\documentclass[10pt,oneside,a4paper,leqno,final]{article}
\usepackage{latexsym}
\usepackage[latin1]{inputenc}
\usepackage{color}
\usepackage[active]{srcltx}
\usepackage{multirow,rotating}
\usepackage{a4wide}
\usepackage[english]{babel}
\usepackage{amscd}
\usepackage{epsfig}
\usepackage{subfigure}
\usepackage{graphicx}
\usepackage{amsmath}
\usepackage{amssymb}
\usepackage{mathrsfs}
\input xy
\xyoption{all}
\usepackage{geometry}

\newtheorem{thm}{\sc Theorem}[section]      
\newtheorem{lem}[thm]{\sc Lemma}            
\newtheorem{prop}[thm]{\sc  Proposition}     
\newtheorem{defn}[thm]{\sc Definition}      
\newtheorem{rem}[thm]{\sc Remark}       

\newcommand{\RR}{\mathbb{R}}
\newcommand{\CC}{\mathbb{C}}

\newcommand{\from}{\colon}
\newcommand{\C}{\mathbb C}

\newcommand{\N}{I\!\!N}
\newcommand{\R}{\mathbb R}
\newcommand{\Z}{\mathbb Z}

\newcommand{\n}[1]{{\bf #1}}

\newcommand{\proof}{{\sl Proof.}\hspace{5pt}}   
\newcommand{\finedim}{\hfill $\Box$}            

\newcommand{\set}[2]{\left\{{#1}\mid{#2}\right\}}       

\usepackage{bm}
\newcommand{\simbolovettore}[1]{{\boldsymbol{#1}}}

\newcommand{\vC}{\simbolovettore{C}}

\newcommand{\ve}{\simbolovettore{e}}

\newcommand{\vh}{\simbolovettore{h}}

\newcommand{\vM}{\simbolovettore{M}}

\newcommand{\vp}{\simbolovettore{p}}
\newcommand{\vP}{\simbolovettore{P}}
\newcommand{\vq}{\simbolovettore{q}}
\newcommand{\vQ}{\simbolovettore{Q}}

\newcommand{\vR}{\simbolovettore{R}}
\newcommand{\vs}{\simbolovettore{s}}
\newcommand{\vS}{\simbolovettore{S}}

\newcommand{\vu}{\simbolovettore{u}}

\newcommand{\vv}{\simbolovettore{v}}

\newcommand{\vz}{\simbolovettore{z}}

\newcommand{\vxi}{\simbolovettore{\xi}}

\newcommand{\zero}{\boldsymbol{0}}

\newcommand{\norm}[1]{\protect\left\protect\Vert\protect#1\protect\right\protect\Vert}

\begin{document}
\pagenumbering{arabic}

\title{%
  Global dynamics of stationary, dihedral, nearly-parallel vortex
  filaments.}

\author{Francesco Paparella\thanks{Dipartimento di Matematica, Università del Salento, Italy.
{\bf email\/}: \texttt{francesco.paparella@unisalento.it}.}, 
Alessandro Portaluri\thanks{Dipartimento di Matematica, Università del Salento, Italy.
{\bf email\/}: \texttt{alessandro.portaluri@unisalento.it}. Work
    partially supported by``Progetto 5 per mille per la ricerca''
    (Bando 2011). ``Collisioni fra vortici puntiformi e fra filamenti
    di vorticita': singolarita', trasporto e caos.''. Work partially supported by 
the PRIN2009 grant "Critical Point Theory and Perturbative Methods for Nonlinear Differential Equations".}}

\date{\today}
\maketitle

\begin{abstract}
  The goal of this paper is to give a detailed analytical description
  of the global dynamics of $N$ points interacting through the
  singular logarithmic potential and subject to the following symmetry
  constraint: at each instant they form an orbit of the dihedral group
  $D_l$ of order $2l$. The main device in order to achieve our results
  is a technique very popular in Celestial Mechanics, usually referred
  to as {\em McGehee transformation.\/} 
  After performing this change
  of coordinates that regularizes the total collision, we study the
  rest-points of the flow, the invariant 
  manifolds and we derive interesting information about the
  global dynamics for $l=2$.  We observe that our problem is
  equivalent to studying the geometry of stationary configurations of
  nearly-parallel vortex filaments in three dimensions in the LIA
  approximation.
\end{abstract}
\noindent {\em MSC Subject Class\/}: Primary 70F10; Secondary 37C80.
\vspace{0.5truecm}

\noindent {\em Keywords\/}: Dihedral $N$-vortex filaments, McGehee
coordinates, global dynamics.

\section*{Introduction} \label{sec:intro}

Equations of motion for interacting point vortices were introduced
by Helmholtz in a seminal paper published in 1858. Towards the end of
the paper he introduced the point vortex model, by considering the
trace of the point of intersection of a family of infinitely thin,
straight parallel vortex filaments with a plane perpendicular to one
(and then to all) of them. A large and still growing body of
literature has focused on the study of point vortices, seen as
vorticity monopoles in a two-dimensional ideal fluid. One may think of
them as playing a role in ideal hydrodynamics similar to that played
by point masses in Celestial Mechanics. For a review see the book by
Newton \cite{Newton} and references therein.

A comparatively much smaller body of literature has been devoted to
the study of vortex filaments in a three-dimensional ideal fluid, a
subject probably closer to Helmholtz's original ideas. Based upon the
equations governing the evolution of vorticity in three dimensions, Da
Rios derived in 1906 the localized self-induction approximation (LIA)
describing the approximate motion of an isolated filament and
later re-derived by Arms and Hama in 1965.  In 1972 Hasimoto introduced a
change of coordinates which takes the familiar Frenet-Serret formulas
for the geometry of a curved filament and the localized self-induction
equation governing its dynamics into the nonlinear Schr\"odinger
equation, a completely integrable infinite dimensional Hamiltonian
system.  Building on these results, Klein, Majda and Damodaran in
\cite{klein} derived in 1995 a simplified model describing the time
evolution of $N$-vortex filaments nearly but not perfectly parallel to
the $z$-axis.  According to this model the motion of $N$-interacting
nearly parallel filaments is given by the following coupled partial
differential equations:
\begin{equation}\label{eq:scroedingernonlineare}
\dfrac{1}{i}\partial_t\Psi_j\,=\,\tilde \Gamma_{j}\partial_{\sigma}^{2}\Psi_{j}+ \sum_{k\neq j}\, \dfrac{\tilde \Gamma_{k}}{2\pi}
\dfrac{\Psi_{j}-\Psi_{k}}{\|\Psi_{j}-\Psi_{k}\|^{2}}
\end{equation}
for $\n{N}:=\left\{ 1,...,N\right\} $, $\tilde \Gamma_{j}\in\mathbb{R^*}$
and where each  $\Psi_{j}$ is a  complex-valued function. 
Among all the solutions of the system \eqref{eq:scroedingernonlineare}
a special role is played by the {\em stationary solutions\/}.

In this paper we are interested in studying the geometry of the
stationary solutions of \eqref{eq:scroedingernonlineare}.  Therefore,
by substituting in \eqref{eq:scroedingernonlineare} $\Psi_{j}(\sigma,
t)$ with $\vq_j(\sigma)$, where $\vq:\mathbb{R}\to\mathbb{R}^2$ we
reduce \eqref{eq:scroedingernonlineare} to the following system of
coupled ordinary differential equations for stationary vortices
filaments:
\begin{equation}\label{eq:2}
\tilde{\Gamma}_{j}\ddot\vq_j(\sigma)+\sum_{k\neq j}\dfrac{\tilde{\Gamma}_{k}}{2\pi}
\dfrac{\vq_{j}(\sigma)-\vq_{k}(\sigma)}{\|\vq_{j}(\sigma)-\vq_{k}(\sigma)\|^{2}}=0
\end{equation}
where $\dot{}$ denotes differentiation with respect to the arc-length
parameter $\sigma$. 
By multiplying each equation in \eqref{eq:2} by $\tilde \Gamma_j$ and
defining the {\em potential function\/}  $U$ as
\begin{equation}\label{eq:statpotfunc}
U(\vq):=- \sum_{\substack{j=1\\k\neq j}}^N \dfrac{\tilde{\Gamma}_j \tilde{\Gamma}_k}{4\pi} \log \|\vq_j-\vq_k\|
\end{equation}
where $\vq:=(\vq_1, \dots, \vq_N)$, the ODE in \eqref{eq:2} can be written as
\begin{equation}\label{eq:Newton}
\Gamma \ddot \vq = \dfrac{\partial U}{\partial \vq}.
\end{equation}
Here $\Gamma$ is the diagonal block matrix defined by $\Gamma
=[\Gamma_{ij}]$ and $\Gamma_{ij}= \Gamma_i^2 \delta_{ij} I_2$ where
$I_2$ denotes the two by two identity matrix.  By interpreting the
parameter $\sigma$ as a time-like coordinate, the equation
\eqref{eq:Newton} can be seen as describing the dynamics of point
masses on the plane, interacting with a logarithmic central potential.
Therefore, from now on, we shall refer to the parameter $\sigma$ as a
{\em time-parameter.\/} Let $\Gamma^*= \sum_{j=1}^N \tilde
\Gamma_j^2$.  We define the {\em center of vorticity\/} $C$ as $ \vC:=
\dfrac{1}{ \Gamma^*}\sum_{j=1}^N \tilde\Gamma_j^2 \vq_j.$ From
\eqref{eq:Newton} it is straightforward to check that the center of
vorticity satisfies the conservation law $\ddot{C}=0$.  A severe
difficulty in investigating a system such as \eqref{eq:Newton} is due
to the presence of singularities both partial and total: physically
they represent collisions between some or all of the vortices. This
problem motivates the introduction of a change of coordinates that
regularizes the equation of motion. A change of coordinates of this
sort was recently used by Stoica \& Font in \cite{stoica} for a
singular central log-potential problem which arises in galactic
astrodynamics. This is an example of {\em McGehee transformations\/},
a regularizing change of variables currently popular in the field of
Celestial Mechanics and first introduced in 1974 by
McGehee\cite{McGehee74} in order to study orbits passing close to the
total collapse in the collinear three-body problem in $\mathbb{R}^3$.

McGehee transformations consists of a polar-type change of coordinates
in the configuration space, a suitable rescaling of the momentum and a
time scaling. The idea behind this non Hamiltonian change of
coordinates is to blown up the total collision to an invariant
manifold called {\em total collision manifold\/} over which the flow
extends smoothly. Furthermore, each hypersurface of constant energy
has this manifold as a boundary.  The effect of rescaling time, is to
study some qualitative properties of the solutions close to total
collision. In fact, by looking at the transformation defined in
equation \eqref{eq:sz}, it readily follows that the effect of this
transformation is to slow the motion in the neighborhood of the total
collapse, which is reached in the new time asymptotically.

However, we point out that there are several crucial differences
between the classical McGehee transformations in use in celestial
mechanics and those appropriate for our case. The most important one
is related to the lack of homogeneity of the logarithmic nonlinearity.
This breaks down some nice and useful properties of the
transformation. For instance it is not possible as in the $N-$body
problem, to recover the global dynamics of orbits passing arbitrarily
close to total collapse by merely looking at a suitable Lyapunov
function defined on the total collision manifolds whose only rest
points represent the central configurations of the bodies in
gravitational interaction. (Compare \cite{McGehee74}, \cite{Devaney80}
\cite{Devaney81}, for further details).

Nevertheless, it is still possible to regularize the vector field and
therefore we can still carry out a detailed analytical description of
the rest points, of the invariant manifolds and of the
heteroclinics. The ability of these coordinates to give insight on the
dynamical properties of the problem becomes apparent when studying a
class of symmetric solutions that we call {\em dihedral equivariant}
(Section \ref{sec:equivsetup}). In the new coordinates we were able to
investigate the global dynamics of the problem and to prove some
important features which were not at all obvious when stating the
problem in conventional Cartesian coordinates. For example, we could
prove that a dihedral equivariant configuration of four vortices is
always bounded (Theorem \ref{thm:nohyp}), even if the standard
argument based on total energy conservation is unable to rule out
unbounded solutions, as it does for central potentials like those
studied in \cite{stoica}.

Another important feature is the fact that every dihedral equivariant
solution experiences a collision within a finite time (Theorem
\ref{thm:sibinary}). This may be either a total collapse (simultaneous
collision of all the vortices) or a binary collision (simultaneous
pairwise collisions of the vortices). Only the second type of
solutions is generic, in the sense that it is generated by a set of
initial conditions having full Lebesgue measure (Theorem
\ref{thm:6A}). By using a recent result proven in \cite{cate} we were
able to study the dynamics for arbitrarily long times, by defining
generalized solutions that continue by {\em transmission\/} after a
binary collision.


On the other hand, an interesting qualitative feature that the present
problem has in common with the gravitational $n$-body problem, is that
the total collapse can only be reached in central configuration. In
Celestial Mechanics this is a well-known fact first proven by Sundman
for the three body problem and generalized by Wintner in the general
case. Recently, the authors of \cite{gencoll} proved that this result
also holds for a class of weak potentials, including the logarithm. We
recover this property as a dynamical feature of the equations written
in McGehee coordinates. It is worth noticing that for genuine point
vortices (and not for stationary configurations of vortex filaments,
as in our case) total collapse need not happen along central
configurations (that is, with self-similar motion; see \cite{Newton}
for further details).


Finally, our results hint at the presence of a chaotic dynamics of the
generalized solution. A challenging problem would be to investigate
the symbolic dynamics or to say something about the topological or
geodesic entropy after a global regularization of the singularities.
The dynamics of the non stationary solutions of the PDEs
\eqref{eq:scroedingernonlineare} also remains to be studied, even for
initial conditions close to the stationary solutions investigated in
the present work.\\
The paper is organized as follows:
\tableofcontents

\section{McGehee coordinates and regularization }

In this section we develop the McGehee--like transformations that we
need in order to study the dynamics close to total collapse. In the
first subsection, we describe the general framework and we fix
notations.  Then, in the following subsection we repeatedly perform
coordinate changes, and we re-parameterize time, until we arrive at a
set of equation of motions which is suitable to study the dynamics of
our problem.

\subsection{General set-up and McGehee coordinates}
Let $N \geq 2$ be an integer.  Let $\zero$ denote the origin in $\R^2$
and let $\tilde \Gamma_1, \dots ,\tilde \Gamma_N $ be $N$ positive
numbers (which can be thought as strength of the vortex filaments).
The conservation law of the center of vorticity implies that there
is always an inertial reference frame where the position of the center
of vorticity is at the origin. In this reference frame we can identify the
configuration space $\vQ$ for the problem with the subspace of $\mathbb R^{2N}$
defined by
\[ 
\vQ:= \{\vq=(\vq_1, \dots, \vq_N) \in \mathbb{R}^{2N}: \ \ \sum_{j=1}^N  \Gamma_j\, 
\vq_j=\zero\}.
\]
For each pair of indexes $i, j \in \n{N}$ let $\Delta_{i,j}$ denote
the collision set of the $i$-th and $j$-th vortex filament; namely
$\Delta_{i,j}= \{\vq \in \mathbb R^{2N}: \vq_i=\vq_j\}$.  We call {\em
  collision set\/} the subset of the configuration space given by
$\Delta:= \bigcup_{i\neq j} \Delta_{i,j}$ and {\em reduced
  configuration space\/}, the set $\hat \vQ:= \vQ\setminus\Delta$. The
Newton's equations in Hamiltonian form can be written as:
\begin{equation}\label{eq:ham}
\left\{
\begin{array}{ll}
\Gamma \dot \vq =\vp\\ &\\
\dot\vp = \dfrac{\partial U}{\partial \vq}
\end{array}\right.
\end{equation}
where the Hamiltonian function $H$ is defined by:
\begin{equation}\label{eq:hamfuncgeneral}
H:T^*(\hat \vQ) \longrightarrow \R \qquad H(\vq, \vp):= \dfrac12 \langle \Gamma^{-1} \vp,
\vp\rangle - U(\vq), 
\end{equation}
on the {\em phase space\/} $T^* (\hat \vQ)$ (the cotangent
bundle over the configuration space). Explicitly $T^*(\hat \vQ)=\hat \vQ \times
\vP$ for 
\[
\vP:= \{\vp=(\vp_1, \dots, \vp_N)\in \mathbb R^{2N}| \ \ \sum_{j=1}^N \vp_j=0\}. 
\]
The differential equations in \eqref{eq:ham} then determine a vector
field with singularities on $\mathbb R^{2N}\times \mathbb R^{2N}$, or
a real analytic vector field without singularities on $(\mathbb
R^{2N}\setminus \Delta)\times \mathbb R^{2N}$.  
The vector field given
by \eqref{eq:ham} is everywhere tangent to $\vQ \times \vP$ and so
this $4(N-1)$ dimensional linear subspace is invariant under the
flow. We henceforth restrict our attention to the flow on the phase
space $\vQ\times \vP$. Consequently $H$ is an integral of the
system. This means that the level sets
$\Sigma_h:=H^{-1}(h)\cap(\vQ\times \vP)$ are also invariant under the
flow \eqref{eq:ham}. We observe that $\Sigma_h$ is a real analytic
submanifold of $\hat \vQ\times \vP$ having dimension $4(N-1)-1$. The
flow, however, is not complete. In fact certain solutions run off in
finite time. This happens in correspondence to any initial condition
leading to a collision between two or more vortex filaments: the
corresponding solution  meet $\Delta$ in finite time.  We shall
call {\em total collapse}, or {\em total collision} the simultaneous
collision of all the vortices. Because the center of vorticity has
been fixed at the origin, if a total collapse happens, it must occur
at the origin of $\vQ$.

We introduce the {\em angular momentum}. Denoting by $(\R^2, \omega_0) $ the two dimensional
symplectic vector space, the angular momentum is given in terms of the
standard symplectic structure as
\begin{equation}\label{eq:angmom}
\vM\,:=\,\sum_{j=1}^N \Gamma_j^{-1} \, \omega_0(\vq_j, \vp_j) .
\end{equation}
We observe that if $\hat \vq_j:=(\vq_j, 0)\in \R^3$, $\hat\vp_j:=(\vp_j,0) \in \R^3$ 
then the angular momentum defined in \eqref{eq:angmom} agrees with 
\begin{equation}\label{eq:angmom2}
\vM \, \ve_3\,=\, \sum_{j=1}^N  \Gamma_j^{-1} \,(\hat \vq_j \times \hat \vp_j).
\end{equation}

An important class of solutions for many problems of celestial
mechanics is called {\em self--similar\/}. This class is characterized
by the fact that up to rotation and homothecy the configuration of all
point masses is constant in time. Self-similar solutions exist also
for our problem. Before proceeding further, let us first define a {\em central configuration.\/}

\begin{defn}\label{def:cc}
  A point $\vs \in \mathbb{R}^{2N}$, such that $\langle\Gamma \vs, \vs\rangle=1$ is
  called {\em central configuration\/} if
\[
\dfrac{\partial U}{\partial \vq}(\vs) = \mu\ \Gamma\ \vs, \qquad
\mathit{for\ some}\ \ \mu \in \R.
\]
\end{defn}
We now characterize the class of self--similar solutions.
\begin{lem}\label{thm:selfsimilar}
  Let $\vs$ be a central configuration, and $\rho$ be a positive
  smooth function such that
\[
 \ddot \rho =\mu \rho^{-1}.
\]
  Then $\vq(\sigma)=\rho(\sigma)\, \vs$ is a solution of
\eqref{eq:Newton}.
\end{lem}
\proof 
By a direct substitution, we have
\[
\Gamma\, \ddot \vq \,=\, \Gamma \,\ddot \rho\, \vs \,= 
\,\Gamma \mu\, \rho^{-1} \vs\,=
\, \rho^{-1} \dfrac{\partial U}{\partial \vq}(\vs)=\dfrac{\partial U}{\partial \vq}(\vq) ,
\]
where the last equality follows by the fact that the gradient of $U$
is a homogeneous function of degree $-1$.\finedim

Let $\varphi_1, \varphi_2$ be two monotonically increasing, unbounded,
smooth functions on the positive half-line, vanishing when $r\to 0^+$.
We define the following McGehee-like coordinates:
\begin{equation}\label{eq:sz}
\left\{
\begin{array}{ll}
\norm{\vq}= \varphi_1(r)\\ \\
\vs=\vq/\norm\vq\\ \\
\vz = \varphi_2(r) \vp.  
\end{array}
\right.
\end{equation}
The vector $\vs$ lies on the $2(N-1)-1$-dimensional unit round sphere
centered at the origin, that we shall call the {\em shape sphere}
$\vS$.  As a direct consequence of the fact that the center of
vorticity is fixed at the origin, it follows that a {\em total
  collapse} may only happen at $r=0$. This constraint defines the {\em
  collision manifold.}

\subsection{Equation of motions in $(r,\vs, \vz)$- coordinates}
In order to write the equations of motion using the McGehee
coordinates, it will be useful to define $v:= \langle \Gamma^{-1} \vz,
\vs \rangle$ and to rescale 
time according to $[\varphi_1\varphi_2]^{-1}d\sigma=\,d\tau$.
In these new coordinates
the submanifold of constant energy $h$ is given by
\begin{equation}\label{eq:stratoenergiacostsz}
  \Sigma_h=\set{(r,\vs, \vz)\in \R^*_+\times \vS \times \R^2}{
    \langle \Gamma^{-1} \vz, \vz \rangle \, =\, 2 \varphi_2^2(r)\left[h+  U(\varphi_1(r)\vs)\right]}.
\end{equation}
By taking into account the expression \eqref{eq:statpotfunc} for the
potential, we have
\[
\begin{split}
\varphi_2^2(r) U(\varphi_1(r)\vs)&:= -\varphi_2^2(r)\sum_{k\neq j} \dfrac{\tilde
{\Gamma}_j \tilde{\Gamma}_k}{4\pi} \log \|\vq_j-\vq_k\|=\\
 &=-\varphi_2^2(r)\left[\sum_{k\neq j} \dfrac{\tilde{\Gamma}_j
\tilde{\Gamma}_k}{4\pi}
\left(\log \varphi_1(r)+ \log \|\vs_j-\vs_k\|\right)\right].
\end{split}
\]
In order $\Sigma_h$ to be not empty we choose the functions $\varphi_j$ as follows:
\[
\varphi_1(r):= \eta(r) e^{-[\varphi_2(r)]^{-2}}, \qquad  \eta(r)\in O(r^k), 
 \; k>0.
\]
With this choice, we get $1:=\lim_{r \to 0^+} - \varphi_2^2(r) \log\varphi_1(r)$
and by this  it follows that the stratum $\Sigma_h$ meets the
collision manifold $r=0$  along the submanifold
\[
\Lambda:=\set{(0,\vs, \vz)\in \R^*_+\times \vS \times \R^2}{\dfrac12 \langle \Gamma ^{-1}
\vz, \vz \rangle = \, G},
\]
where $G:=\sum_{k< j} \dfrac{\tilde{\Gamma}_j
\tilde{\Gamma}_k}{2\pi}$.  We observe that $\Lambda$ does not depend
on the energy level $h$. Therefore, all strata $\Sigma_h$ share the
same boundary at the collision manifold.
Following the authors in \cite{stoica}, one possible choice is the following:
\begin{equation}\label{eq:lefi}
\left\{
\begin{array}{ll}
\varphi_1(r):= r e^{-1/r^2}\\
\\
\varphi_2(r):= r
\end{array}\right.
\end{equation}
With this choice, the Hamiltonian system (\ref{eq:Newton}) becomes the
following system of ordinary differential equations
\begin{equation}\label{eq:mcgehee1bisdie2}
\left\{
\begin{array}{ll}
 r' = \dfrac{r^3}{r^2+2}\, v\\
\\
\vs' =\,\Gamma^{-1}\,\vz- v \vs\\
\\
\vz' = r^2\left[\dfrac{1}{r^2+2} v \vz+\dfrac{\partial U}{\partial \vq}(\vs)\right].\\
\end{array}\right.
\end{equation}
where, in the last equation, we have used the identity
$re^{-1/r^2}\partial_\vq U(\vq)= \partial_\vq U(\vs)$ which follows
from the homogeneity of degree $-1$ of $\partial_\vq U$, together with
(\ref{eq:sz}) and (\ref{eq:lefi}).

\section{The dihedral problem}
\label{sec:equivsetup}

\begin{figure}
\centering
\subfigure[$l=2$]{%
\includegraphics[width=0.480\textwidth]{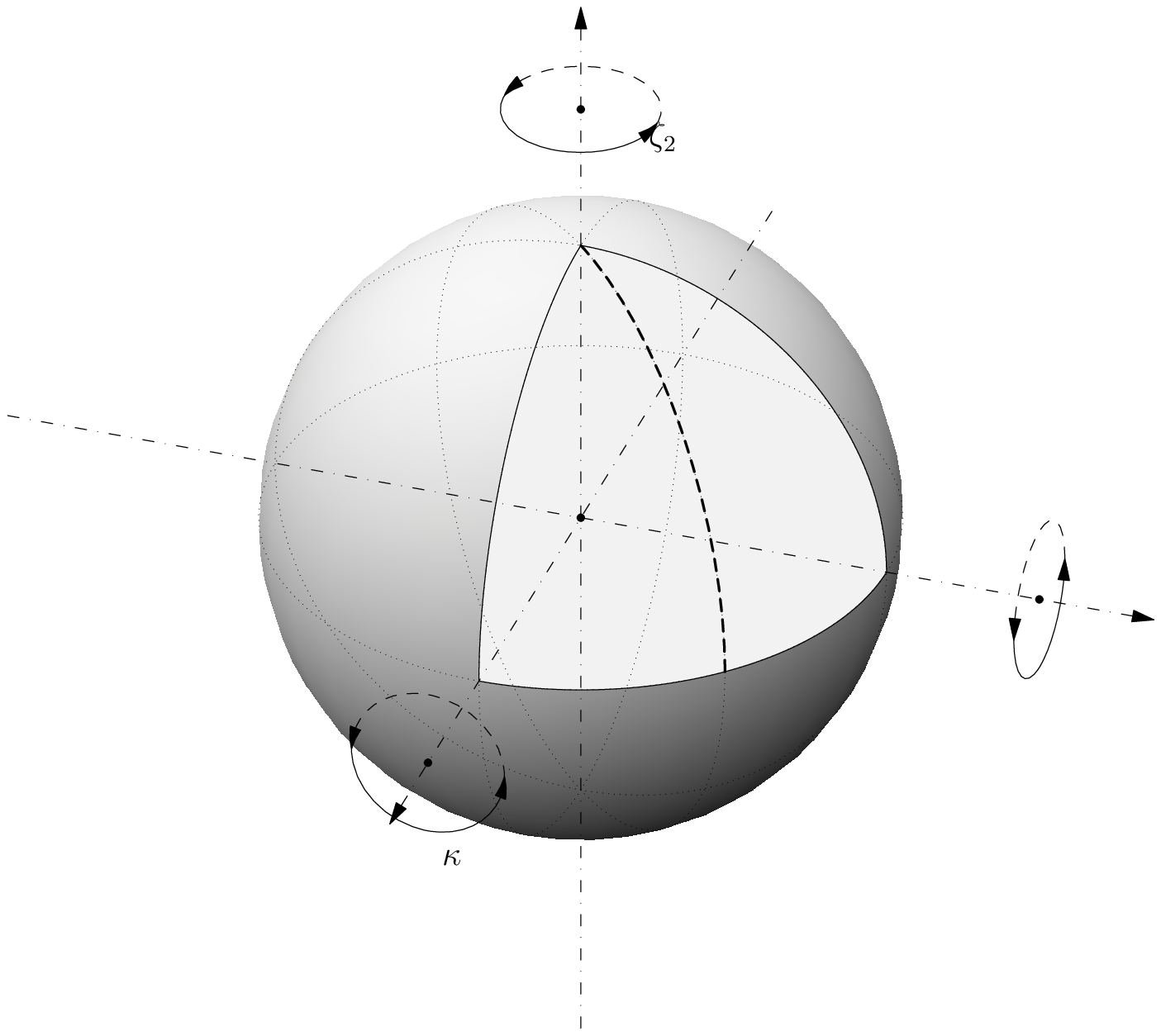}
\label{fig:fund4}
}
\subfigure[$l=3$]{%
\includegraphics[width=0.480\textwidth]{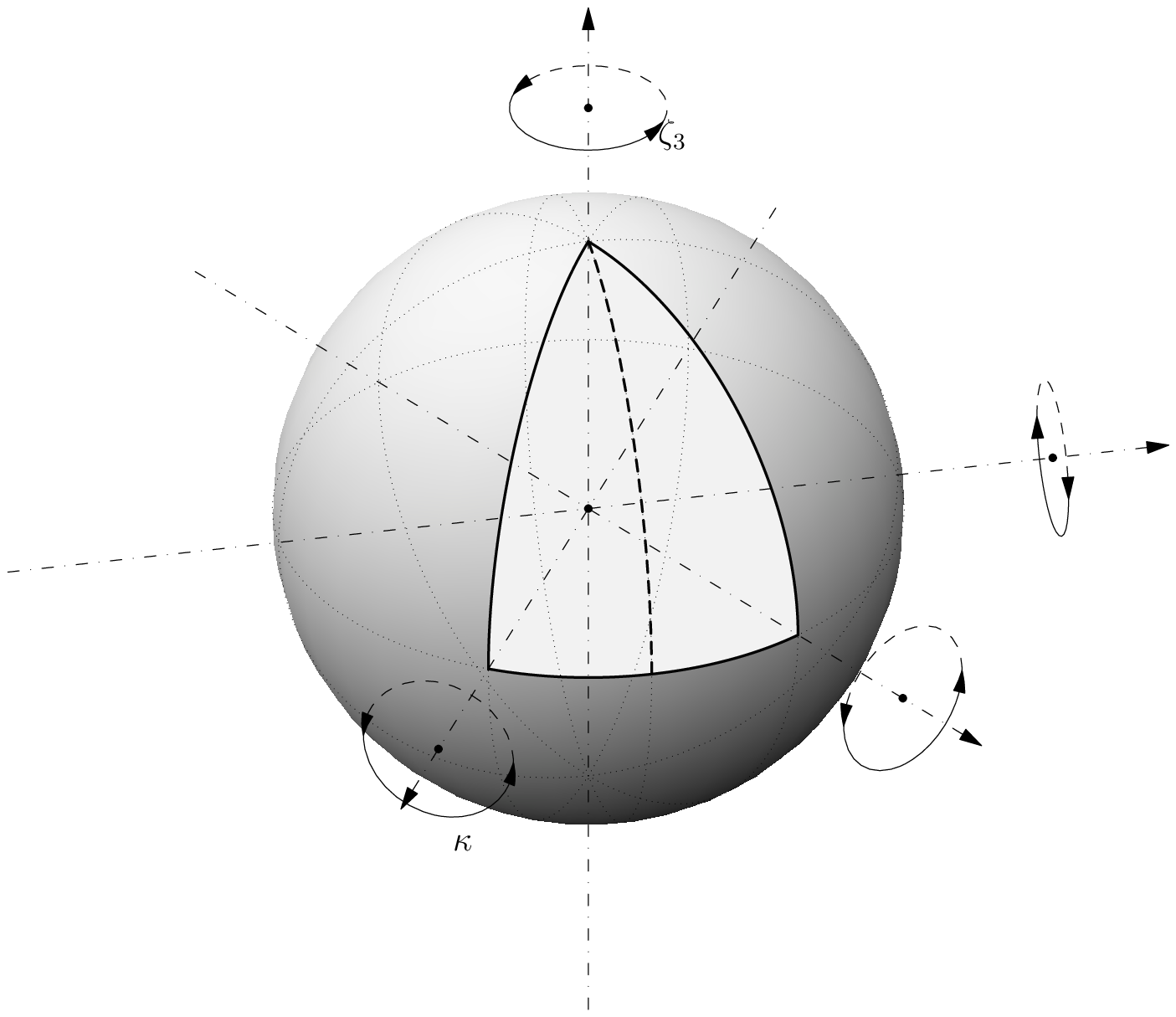}
\label{fig:fund}
}
\caption{Dihedral groups $D_l$, with the upper half of the fundamental
domains in white.}
\label{fig:fun}
\end{figure}

Let us recall some basic facts about the
dihedral group, seen as a map of $\mathbb R^3$ into itself. For
further details we shall refer to \cite{fp}.  Let $\R^3 \cong \CC
\times \RR$ be endowed with 
coordinates $(z,y)$, $z\in \CC$, $y\in
\RR$. For $l\geq 1$, let $\zeta_l$ denote the primitive root of 
unity
$\zeta_l = e^{2\pi i/l}$; the \emph{dihedral} group $D_{l}\subset
SO(3)$ is the group of order $2l$ 
generated by 
\[\label{eq:dihedral}
\zeta_l\from (z,y) \mapsto (\zeta_l z, y ),\qquad
\kappa\from (z,y) \mapsto (\overline{z},-y),
\]
where $\overline{z}$ is the complex conjugate of $z$. The
non-trivial elements of $D_l = \langle \zeta_l,\kappa \rangle$ are
the  $l-1$ rotations around the $l$-gonal axis $\zeta_l^j$,
$j=1,\ldots,l-1$ and  the $l$ rotations of angle $\pi$ around the
$l$ digonal  axes orthogonal to the $l$-gonal axis (see figure
\ref{fig:fun}) $\zeta_l^j\kappa$, $j=1,\ldots ,l$. In
figures~\ref{fig:fund4} and~\ref{fig:fund} we show the
upper-halves of the fundamental domains for the action of $D_l$
restricted on the unit sphere for $l=2$ and $l=3$, respectively. 

The action of $D_l$, restricted on the fixed subspace
$\left(\R^{4l}\right)^{D_{l}} \cong \R^2$, generates a dihedral
configuration of vortices for any given $\vq_0\in \R^2$. If we assign
the same circulation to all vortices, then we can express the
logarithmic potential acting on a single vortex as 
\begin{equation}\label{eq:potenziale}
U(\vq_0) = -\sum_{g\in D_l \smallsetminus \{1\}}
\log(\left\| \vq_0 - g \vq_0 \right\|),
\end{equation}
where, without further loss of generality, we have taken $\tilde
\Gamma_i^{2}=2\pi/l$.  
A solution of equations \eqref{eq:ham} such
that the vortices lie on a dihedral configuration at all times, is a
{\em dihedral equivariant orbit}.

\begin{lem}\label{thm:angmomiszero}
The angular momentum of any dihedral equivariant orbit is zero. 
\end{lem}
\proof Using \eqref{eq:angmom2} and Hamilton's equations \eqref{eq:ham},
it follows that the angular momentum is:
\[
\vM\,= \,\sum_{j=1}^N\, \hat \vq_j \times \dot{\hat\vq}_j.
\]
By hypothesis, at each instant the configuration of the vortices is
given by the action of the dihedral group $D_l$; $\vq_j= g_j \vq$ for
$g_j \in D_l$. Hence
\[
\hat \vq_j \times \dot{\hat\vq}_j= g_j\hat  \vq \times g_j\dot{\hat\vq}\,=\, g_j\big(\hat  \vq \times \dot{\hat\vq}\big).
\]
We observe that $\hat \vq \times \dot{\hat \vq} = (0,0,y)$, for some
real number $y$. Therefore, using (\ref{eq:dihedral}), we have
\[
\vM=\sum_{j=1}^{2l} g_j \big(\hat  \vq \times \dot{\hat\vq}\big) = 
l(0,0,y) + l(0,0,-y) = \simbolovettore{0}.
\]
This conclude the proof.\finedim
\begin{rem}
  As consequence of the Sundman-type estimates proved by authors in
  \cite{gencoll}, for a large class of potentials including the
  logarithmic one, a necessary condition in order to have the total
  collapse is that the total angular momentum of the system should be
  zero. Therefore as direct consequence of Lemma
  \ref{thm:angmomiszero}, we can conclude that total collision orbits
  can occur.
\end{rem}
\subsection{The logarithmic dihedral potential}
For $z=\rho \,e^{i\alpha}$, with $\rho\ge0$ and $\alpha_j=2\pi j/l$, $j=1,\ldots,l-1$, it is easy
to verify that:
\[
\begin{aligned}
|z-\bar z|&= 2 \rho \left|\sin \alpha\right|\\
|z-\zeta_l^jz|&= 2 \rho \left|\sin\left(\dfrac{\alpha_j}{2}\right)\right|\\
|z-\zeta_l^j\bar z|&= 2 \rho \left|\sin\left(\dfrac{\alpha_j-2 \alpha}{2}\right)\right|.\\
\end{aligned}
\]
Therefore, for $\alpha \in (0, \pi /l)$ (that is, for $\vq$ in the
fundamental domain) the potential \eqref{eq:potenziale} becomes:
\begin{equation}\label{eq:logpot}
U(\rho, \alpha):= -\log\left((2\rho)^{2l-1}\right)-\log \left[\sin\alpha\, \prod_{j=1}^{l-1} \sin\left(\dfrac{\alpha_j}{2}\right)\,
\sin\left(\dfrac{\alpha_j}{2}-\alpha\right)\right].
\end{equation}
Since $\prod_{j=1}^{l-1}\sin\left(j\dfrac{\pi}{l}\right)=\dfrac{l}{2^{l-1}}$ and by taking into account the multiple-angle 
formula
\[
 \sin(n x) = 2^{n-1} \prod_{k=0}^{n-1} \sin\left(\dfrac{k\pi}{n}+x\right),
\]
we have 
\[
 \sin\alpha\prod_{j=1}^{l-1}\sin\left(\dfrac{j\pi}{l}-\alpha\right)= \dfrac{\sin(l\alpha)}{2^{l-1}}.
\]
These calculations lead to the following definition.
\begin{defn}\label{def:potenziali}
We define the {\em  dihedral logarithmic potential\/} as :
\begin{equation}\label{eq:logpotdiedralunred}
U(\vq):= -\log \rho^{2l-1}-\log \left(2l\sin(l\,\alpha)\right),
\end{equation}
and its angular part as:
\begin{equation}\label{eq:logpotdiedral}
U(\alpha):= -\log \left(2l\sin(l\,\alpha)\right).
\end{equation}
\end{defn}
\begin{lem}\label{thm:plcc} (Planar $2l$-gon) The dihedral problem admits 
exactly one (up to permutation of the vortices) central
configuration, which is given by the vertices $(e^{(2k+1)\pi i/(2l)},0)$
of a regular $2l$-gon.
\end{lem}
\proof From the definition (\ref{def:cc}) it follows that central
configuration corresponds to a critical point of the potential
restricted to the shape sphere. From (\ref{eq:logpotdiedral})
it follows that
\[
 U'(\alpha)=-\dfrac{l\,\cos(l\alpha)}{\sin(l\alpha)}
\]
whose critical point in the fundamental domain of the shape sphere is
$\alpha_c=\dfrac{\pi}{2l}$. From $U''(\alpha_c)=l^2$
it also follows that the point  $\alpha_c$ is a minimum. \finedim\\
We also observe that
\begin{enumerate}
\item $U(\alpha_c)<0$ for $l=2,3$;
\item  $U(\alpha_c)>0$ for $l \geq 4$.
\end{enumerate}
By using the identity
\[
\sin(n\alpha)=\sum_{k=0}^n \binom{n}{k} \cos^k\alpha \sin^{n-k}\alpha \sin\left(\dfrac12(n-k)\pi\right),
\]
the expression \eqref{eq:logpotdiedral} may be written in terms of the local 
parameterization of the unit sphere $\vs(\alpha)=(s_1, s_2):=(\cos\alpha, \sin \alpha)$ as follows
\begin{equation}\label{eq:potdieins}
U(\vs)=-\log\left[2l\sum_{j=0}^l\binom{l}{j} s_1^j s_2^{l-j} \sin\left(\dfrac\pi2(l-j)\right) \right]. 
\end{equation}
For $j \in \n{2}$,  $q_j=s_j\|\vq\|$ and hence the dihedral potential given in definition  
\ref{def:potenziali}  becomes:
\begin{equation}\label{eq:potdieinq}
U(\vq)=-\log\left[2l\sum_{j=0}^l\binom{l}{j}\|\vq\|^{l-1} q_1^j q_2^{l-j} \sin\left(\dfrac\pi2(l-j)\right) \right]. 
\end{equation}
\subsection{The geometry of the energy hypersurfaces}
From \eqref{eq:stratoenergiacostsz} let us define $\hat E(h, r, \vs):=
2r^2(h+ U(re^{-1/r^2}\vs))$. 
More explicitly $\hat E$ is given by:
\[
\hat E(h,r,\vs):= E(h,r)+ 2r^2\,U(\vs),
\]
where  $E$ is given by
\begin{equation}\label{eq:energyradial}
E(h, r):= 2\left[h\, r^2 +(2l-1)(1-r^2\log r)\right].
\end{equation}
\begin{rem}\label{rem:Ehat_ge_0}
The boundary of the regions where the motion occurs is given by 
$\hat E(h,r,\alpha)=0$; more precisely, since the kinetic term
$\langle \Gamma^{-1} \vz, \vz \rangle$ is a positive definite
quadratic form, for any fixed energy level $h$ the motion is possible where:
\[
\hat E(h,r,\alpha)\,:=\,2\left[h\, r^2 +(2l-1)(1-r^2\log r)\right]-2r^2\log(2l \sin(l\alpha)) \geq\, 0.
\]
The shape of the curve  $\hat E(h,r,\alpha)= 0$ is key in order to understand
the dynamics of our problem.
\begin{figure}[t]
\centering
{%
\includegraphics[width=0.440\textwidth]{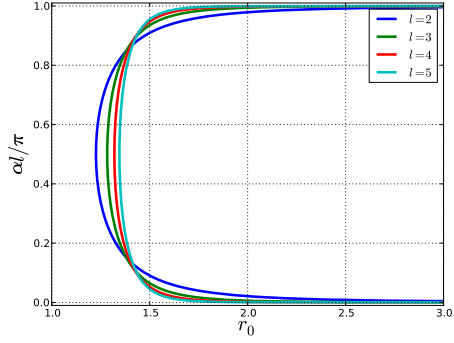}
}
\caption{Zero set of the function $(r,\alpha)\mapsto\hat E(0,r,\alpha)$ and different value of the total number of point vortex.}
\label{fig:plotenrgiazerodiversil}
\end{figure}
\end{rem}
Therefore the hypersurface corresponding to the energy level $h$ is given by
\begin{equation}\label{eq:stratoenergiacostszdie}
\Sigma_h=\set{(r,\alpha, \vz)\in \R^*_+\times (0,\pi/l) \times \R^2}{
\langle \Gamma^{-1} \vz, \vz \rangle \, =\,\hat{E}(h,r,\alpha) }.
\end{equation}
We also observe that $\Sigma_h$ meets the boundary $r=0$ along a submanifold given by
\[
\Lambda:=\set{(0,\alpha, \vz)\in \R^*_+\times (0,\pi/l) \times \R^2}{
\langle \Gamma^{-1} \vz, \vz \rangle \, =\,2\,(2l-1) }.
\]
The set $\Lambda$ is diffeomorphic to an open cylinder and it
represents the component over the fundamental domain of the {\em total
  collision manifold\/}. As we shall see in the following, the total
collision manifold is homeomorphic to a two dimensional torus.
\begin{figure}[ht]
\centering
{%
\includegraphics[width=0.440\textwidth]{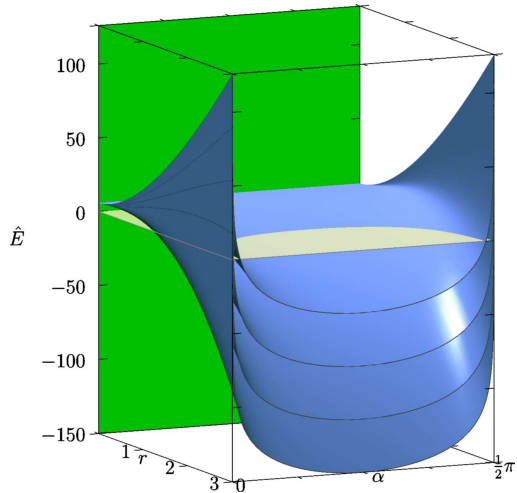}
}
\caption{Energy hypersurfaces for different values of the energy level $h$ and the zero plane.}
\label{fig:energy}
\end{figure}

Let us introduce a further change of coordinates.  In the
dihedral problem the shape sphere reduces to $S^1$, therefore it is
natural to use the parameterization  $\alpha \mapsto \vs(\alpha)$, where 
$\vs(\alpha) = (\cos \alpha, \sin \alpha )$.
Following \cite{stoica}, we exploit the constraint 
(\ref{eq:stratoenergiacostszdie}) to parameterize the
 we also introduce a new local parameterization the momentum 
coordinates $\vz$ with an angle $\psi$ in the following way
\begin{equation}\label{eq:angolozdie}
\vz = \sqrt{\hat E(h, r,\alpha)}\, (\sqrt{\frac{2\pi}{l}}\cos \psi,\sqrt{\frac{2\pi}{l}} \sin \psi).
\end{equation}
Note that with the choice of the strength of the circulations made
in Section \ref{sec:equivsetup}, we have $\Gamma^{-1}= l/2\pi I$. 
We also use a new time-like variable $\zeta$, defined by 
\[d\tau=\sqrt{\hat  E}\,d\zeta.
\label{eq:timescaling}
\]
With the above parameterizations, the system (\ref{eq:mcgehee1bisdie2}) 
reduces to
\begin{equation}\label{eq:mcgehee22bisdie}
\left\{
\begin{array}{ll}
\dfrac{dr}{d\zeta}= \sqrt{\dfrac{l}{2\pi}}\dfrac{r^3}{r^2+2}\, \hat E(h,r, \alpha)\,(\cos\alpha\cos\psi+\sin\alpha\sin\psi)\\
\\
\dfrac{d\alpha}{d\zeta}=\sqrt{\dfrac{l}{2\pi}}\hat E(h,r,\alpha)(\sin \psi \cos \alpha- \sin \alpha \cos \psi)\\
\\
\dfrac{d\psi}{d \zeta} =\sqrt{\dfrac{l}{2\pi}} r^2(\partial_{q_2}U(\vs(\alpha))\,\cos \psi -\partial_{q_1} U(\vs(\alpha))\sin\psi).\\
\end{array}\right.
\end{equation}
The following result is an obvious consequence of equations 
(\ref{eq:mcgehee22bisdie}).
\begin{lem}
The collision manifold at $r=0$ and the zero velocity manifold $\hat E(h,r, \alpha)=0$ 
are invariant manifolds. Moreover:
\begin{enumerate}
\item On the collision manifold $r=0$ the dynamics is given by
\begin{equation*}
\dfrac{d\alpha}{d\zeta} =-\sqrt{\dfrac{l}{2\pi}}2(2l-1)\sin (\psi-\alpha), 
\qquad \dfrac{d\psi}{d\zeta}=0;
\end{equation*}
\item On the zero-velocity manifold $\hat E(h,r,\alpha)=0$, the
  dynamics is obtained by integrating the last equation in
  \eqref{eq:mcgehee22bisdie}, keeping $r,$ $\alpha$ as constants.
\end{enumerate}
\end{lem}

\section{Flow and  invariant manifolds: the Klein group $D_2$}\label{sec:stabunstab}

For $l=2$ (the Klein group) we may carry out a detailed study of the
global dynamic of the problem. We observe that the dihedral
potential reduces to:
\[
U(\vq) =-\log(8 \|\vq\| q_1q_2)
\]
whence it follows $ \partial_{\vq} U(\vs) = -\left(\dfrac{1}{s_1}+ s_1,   \dfrac{1}{s_2}+ s_2\right).$
In this case the equations of motion become:
\begin{equation}\label{eq:mcgehee22trisdiel=2}
\left\{
\begin{array}{ll}
 \dfrac{dr}{d\zeta}= \dfrac{1}{\sqrt \pi} \dfrac{r^3}{r^2+2}\,\hat E(h,r,\alpha)\,\cos(\psi-\alpha)\\
\\
\dfrac{d\alpha}{d\zeta}= \dfrac{1}{\sqrt \pi}\hat E(h,r,\alpha)\, \sin(\psi-\alpha)\\
\\
\dfrac{d\psi}{d \zeta} = \dfrac{- r^2}{\sqrt\pi}\left[\dfrac{2\cos(\psi+\alpha)}{\sin 2\alpha}- \sin(\psi-\alpha)\right].
\end{array}\right.
\end{equation}
The rest points of \eqref{eq:mcgehee22trisdiel=2} correspond to the
solutions of the following systems:
\[
\left\{
\begin{array}{ll}
 r=0\\
\sin(\psi-\alpha)=0
\end{array}\right., \qquad 
\left\{
\begin{array}{ll}
 r=0\\
\hat E=0
\end{array}\right., \qquad 
\left\{
\begin{array}{ll}
\hat E=0\\
\dfrac{2\cos(\psi+\alpha)}{\sin 2\alpha}= \sin(\psi-\alpha)
\end{array}\right..
\]
However, it is readily seen that the second system has no solutions,
as the conditions $\hat E=0 $ and $r=0$ are incompatible. With straightforward 
calculations we obtain the following result.
\begin{lem}\label{thm:restpoints}
  The equilibria of the vector field given in
  \eqref{eq:mcgehee22trisdiel=2} consists of four curves, two
  belonging to the collision manifold and the other two on the
  zero-velocity manifold. In the coordinates $(r, \alpha, \psi)$,
  the curves on the collision manifold are given by
\begin{enumerate}
\item[(i)]  $\mathscr P_1\equiv(0,\alpha,\alpha)$;
\item[(ii)]   $\mathscr P_2 \equiv (0,\alpha, \alpha+\pi)$
\end{enumerate}
The curves on the zero-velocity manifold are given by
\begin{enumerate}
\item[(iii)]  $\mathscr P_3\equiv\big(r,\alpha,\arctan m(\alpha)\big)$;
\item[(iv)]  $\mathscr P_4 \equiv \big(r ,\alpha , \arctan m(\alpha)+\pi\big)$,
where
\[
 m(\alpha):= \dfrac{\cos\alpha(\sin^2\alpha+1)}{\sin\alpha(\cos^2\alpha +1)}.
 \]
 and the pair $(r,\alpha)$ satisfies the equation
 $\hat{E}(h,r,\alpha)=0$, which explicitly reads
\[
\log (4 \sin 2\alpha) = h + 3(r^{-2}-\log r).
\]

\end{enumerate}
\end{lem}

\begin{figure}
\centering
{%
\includegraphics[width=0.400\textwidth]{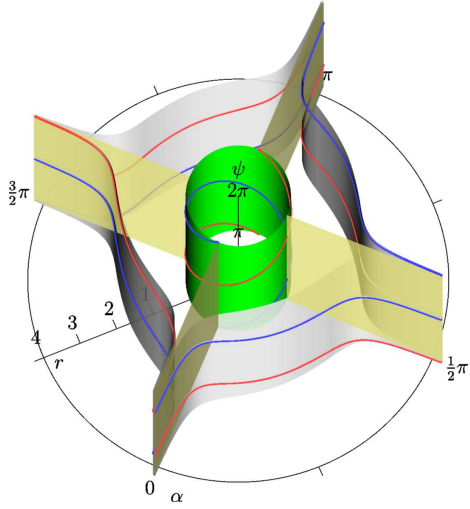}
}
\caption{Rest points and invariant manifolds for the $D_2$ problem. 
The green cylinder is the total collision manifold $r=0$; the blue and 
the red lines lying on it are $\mathscr{P}_1$  and $\mathscr{P}_2$, 
respectively. The four light grey manifolds are the surfaces defined by 
$\hat{E}(0,r,\alpha)=0$ (there is one for each of the four domains of 
the group $D_2$); the red and the blue lines lying on them are 
$\mathscr{P}_3$ and $\mathscr{P}_4$ in the fundamental domain, 
and their images in the other domains.
The blue/red color denotes families of rest points which are transversally 
linearly stable/unstable on their own invariant surface. See the text
for the stability in directions transverse to the invariant surfaces.
The four yellow planes correspond to binary collisions.
\label{fig:zerosetp3eps}}
\end{figure}

The spectrum and the eigenvectors of the Jacobian matrix of
\eqref{eq:mcgehee22trisdiel=2}, evaluated at a rest point, gives
useful information on the local dynamics close to the curves
$\mathscr{P}_1$, $\mathscr{P}_2$, $\mathscr{P}_3$, $\mathscr{P}_4$.
The logic of the calculation is elementary, but the computations
occasionally become too large to be easily manageable. With the help
of a computer algebra system (Maxima 5.21.1) we arrive at the
following results.

\begin{lem}\label{thm:spectrum1}
  The curves of equilibria $\mathscr P_1$ and $\mathscr P_2$ (on the
  collision manifold) given in Lemma \ref{thm:restpoints} are
  degenerate. More precisely:
\begin{enumerate}
\item at each point of the curve $\mathscr P_1$ the spectrum is given
  by $ \mathfrak s(\mathscr P_1)=\{\lambda_1, \lambda_2\}$ where
  $\lambda_1=6/\sqrt \pi$ and $\lambda_2=0$. Furthermore the algebraic
  multiplicity of $\lambda_1$ is $h(\lambda_1)=1$ and the algebraic
  multiplicity of $\lambda_2$ is $h(\lambda_2)=2$. The eigenvectors
  are
\[
 \vv_1=(0,1,0)\ \ \vv_2=(1,0,0)\ \ \vv_3=(0,1,1)
\]
where $\vv_1$ is associated to $\lambda_1$ and $\vv_2$, $\vv_3$ are
associated to $\lambda_2$. The vector $\vv_3$ is tangent to $\mathscr{P}_1$.
\item At each point of the curve $\mathscr P_2$ the spectrum is given
  by $ \mathfrak s(\mathscr P_2)=\{\lambda_1, \lambda_2\}$ where $\lambda_1=-6/\sqrt \pi$ and $\lambda_2=0$.
  The multiplicities and the eigenvectors are the same as 
  for $\mathscr{P}_1$.
\end{enumerate}
\end{lem}
\begin{lem}\label{thm:spectrum2}
  The curves of equilibria $\mathscr P_3$ and $\mathscr P_4$ (on the
  zero velocity manifold) given in Lemma \ref{thm:restpoints} are
  degenerate. More precisely:
\begin{enumerate}
\item at each point of the curve $\mathscr P_3$ the spectrum is given
  by three distinct simple eigenvalues, namely $ \mathfrak s(\mathscr
  P_3)=\{\mu_1, \mu_2, \mu_3\}$ where
\[
\begin{split}
\mu_1&=0\\
\mu_2&=-\dfrac{2r^2\,\big[3\,\cos(\psi-\alpha)e^{6/r^2+2h}+64\,\cos(2\alpha)\,\cos(\psi+\alpha) r^6\big]}{e^{6/r^2+2h}\,\sqrt\pi} \\
\mu_3&= \dfrac{r^2\big[\cos(\psi-\alpha)\,e^{3/r^2+h} + 8 \sin(\psi+\alpha)\,r^3\big]}{e^{3/r^2+h}\,\sqrt\pi},
\end{split}
\]
and 
\[
\hat E_0(r,\alpha)=0, \qquad  \dfrac{2\cos(\psi+\alpha)}{\sin 2\alpha}= \sin(\psi-\alpha).
\]
The eigenvector associated to $\mu_3$ is $\vu_3=(0,0,1)$, which lies
on the invariant manifold $\hat{E}(h,r,\alpha)=0$. The explicit
expressions of the eigenvectors $\vu_1$ and $\vu_2$ associated to
$\mu_1$ and to $\mu_2$ are not short, and we omit them. However,
$\vu_1$ is tangent to the curve $\mathscr{P}_3$ and $\vu_2$ is
transverse to $\hat{E}(h,r,\alpha)=0$.
\item At each point of the curve $\mathscr P_4$ the spectrum is $ \mathfrak s(\mathscr P_3)=\{\mu_1, -\mu_2, -\mu_3\}$. The 
eigenvectors are the same as for $\mathscr{P}_3$.
\end{enumerate}
\end{lem}
\begin{rem}\label{rem:eigvalsign}
  For any $\alpha\in(0,\pi/2)$ it is $\mu_2<0$ and $\mu_3>0$. In fact,
  from the explicit expression of $\psi(\alpha)$ given in
  (\ref{thm:restpoints}) we have $\cos(\psi-\alpha)>0$,
  $\sin(\psi+\alpha)>0$, $\cos(2\alpha)\cos(\psi+\alpha)>0$. The last
  inequality stems from
  \[
  \begin{split}
    \pi/4<\psi(\alpha)<\pi/2-\alpha \quad \mathrm{if} \quad  \hphantom{\pi/}0<\alpha<\pi/4;\\
    \pi/4>\psi(\alpha)>\pi/2-\alpha \quad \mathrm{if} \quad \pi/4>\alpha>\pi/2.
  \end{split}
  \]
\end{rem}
\begin{rem}
  The above results show that the stability properties of the
  equilibria do not depend on the value of the total energy $h$. In
  other words, changing the parameter $h$ does not lead to
  bifurcations. 
\end{rem}

At each rest point let us denote by $W^s$, $W^{u}$ respectively the
(linearly) stable and unstable manifolds, and by $W^0$ the center
manifold. The following results hold.

\begin{prop}\label{thm:dimensionisup3}
  Any rest point in $\mathscr P_3$ and $\mathscr P_4$ is a degenerate saddle. More precisely we have
\begin{enumerate}
\item $
 \dim W^u(\mathscr P_3)=1, \qquad  \dim W^s(\mathscr P_3)=1, \qquad  \dim W^0(\mathscr P_3)=1.$
\item $
\dim W^u(\mathscr P_4)=1, \qquad  \dim W^s(\mathscr P_4)=1, \qquad  \dim W^0(\mathscr P_4)=1.$
\end{enumerate}
\end{prop}
\proof The statement is an immediate consequence of Lemma
\ref{thm:spectrum1}. Note that $\mathscr{P}_3$ and $\mathscr{P}_4$ are
the center manifold of each of their own points, and that they are
{\em neutral}, in the sense that each of their points is an equilibrium.
\finedim\\
\begin{table}\centering
\begin{tabular}{|l|c|c|c|}
\hline & 
\begin{sideways}$\dim  W^{s}\ $\end{sideways} &
\begin{sideways}$\dim  W^{u}\ $\end{sideways} &
\begin{sideways}$\dim W^0\  $\end{sideways}
\\ \hline 
&&&\\
\multirow{1}{*}{At $\mathscr P_1$}&$1$ &--   &$2$  \\ 
\hline  &&&\\
\multirow{1}{*}{At $\mathscr P_2$}&--  &$1$  &$2$  \\ 
\hline  &&&\\
\multirow{1}{*}{At $\mathscr P_3$}&$1$ &$1$  &$1$  \\ 
\hline  &&&\\
\multirow{1}{*}{At $\mathscr P_4$}&$1$ &$1$  &$1$  \\ 
\hline
\end{tabular}
\caption{Dimensions of the invariant  manifolds for the equilibria
  belonging to the curves $\mathscr P_1$, $\mathscr P_2$, $\mathscr P_3$, $\mathscr P_4$}
\label{tb:tabella1}
\end{table}
\begin{prop}\label{pro:piep2selledegeneri}
  For each $\alpha$, the two equilibria
\[
(0,\alpha,\alpha) \in \mathscr P_1
\]
and
\[
(0, \alpha, \alpha+\pi) \in \mathscr P_2 
\]
are nonlinear degenerate saddles. 
\begin{enumerate}
\item $
 \dim W^s(\mathscr P_1)\hphantom{^c}=1, \qquad  \dim W^0(\mathscr P_1)=2.$
\item $
\dim W^u(\mathscr P_2)\hphantom{^c}=1,  \qquad  \dim W^0(\mathscr P_2)=2.$
\end{enumerate}
\end{prop}
\proof The equations of motion for an initial condition on the
collision manifold reduce to
\begin{equation}\label{eq:mcgehee22trisdiel=2collmanfld}
\left\{
\begin{array}{lll}
\dfrac{dr}{d\zeta}=0\\
\\
\dfrac{d\alpha}{d\zeta}= \dfrac{6}{\sqrt \pi}\, \sin (\psi-\alpha)\\
\\
\dfrac{d\psi}{d \zeta} = 0.\\
\end{array}\right.
\end{equation}
It is then evident that the unstable manifold of the equilibrium
$(0,\alpha_0,\alpha_0)\in\mathscr{P}_1$ is the line
$(0,\alpha,\alpha_0)$, and analogously $(0,\alpha, \alpha_0+\pi)$ is
the stable manifold for any point
$(0,\alpha_0,\alpha_0+\pi)\in\mathscr{P}_2$.

The linear stability analysis of Lemma \ref{thm:spectrum1} shows that
each point on $\mathscr{P}_1$ and $\mathscr{P}_2$ has a
two-dimensional center subspace, and hence a two-dimensional center
manifold.  In fact $\mathscr{P}_1$ and $\mathscr{P}_2$ are {\em
  neutral} submanifolds of the center manifold, in the sense that each
of its points is an equilibrium.  

For brevity, here we shall not perform a formal center manifold
reduction. In order to determine the stability of the dynamics on the
center manifold we simply observe that the dynamics projected along
the direction singled out by the eigenvalue $\vv_2=(1,0,0)$, which is
transverse to the total collision manifold, is given by the first
equation in (\ref{eq:mcgehee22trisdiel=2}). From that, we have
$dr/d\zeta>0$ in a neighborhood of any point of $\mathscr{P}_1$ and 
$dr/d\zeta<0$ in a neighborhood of any point of $\mathscr{P}_2$.  

Therefore the equilibria of both $\mathscr{P}_1$ and $\mathscr{P}_2$
are nonlinear saddles.
\finedim


\section{Heteroclinic connections and homothetic orbits}

The existence of heteroclinic connections on and between the invariant
manifolds helps us to develop a global understanding of the flow.

\begin{lem}\label{thm:flussodegcollmnfld}
The flow on the total collision manifold is {\em totally degenerate}. That
is:
\begin{enumerate}
\item[(i)] $W^u(\vP_1) \equiv W^s(\vP_2)$; 
\item[(ii)] $W^u(\vP_2) \equiv W^s(\vP_1)$;
\end{enumerate}
where $\vP_1 \in \mathscr P_1$ and $\vP_2 \in \mathscr P_2$ are chosen
in such a way that the second coordinate of the two points is the
same.
\end{lem}
\proof The proof of this result follows by a straightforward
integration of the equations of motion
(\ref{eq:mcgehee22trisdiel=2collmanfld}) valid on the total collision
manifold $r=0$.  Note that the singularity of the potential at
$\alpha=j\pi/4$, $j=0,\cdots,3$ vanishes on the total collision
manifold, therefore on the heteroclinic connection between $\vP_1$
and $\vP_2$, we have $\alpha\in[0,2\pi)$.
\finedim
\begin{rem}
  By taking into account the identification of the opposite edges of
  the rectangle $(\alpha, \psi)$ with the same orientation, the curves of equilibria
  can be identified with two closed curves on the torus and the
  heteroclinics are represented by arcs joining one point on a
  unstable (red) closed curve with the corresponding point on the
  stable (blue) curve as shown in Figure (\ref{fig:totcollmnfld}).
\begin{figure}
\centering
{%
\includegraphics[width=0.440\textwidth]{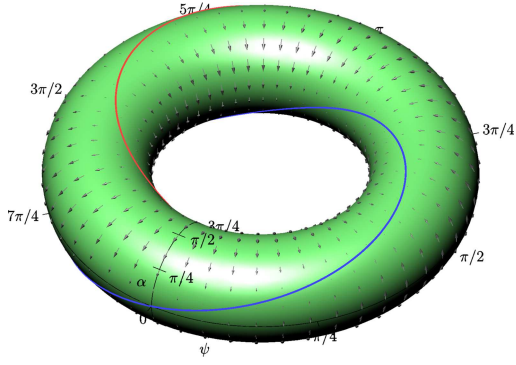}
}
\caption{Total collision manifold and equilibrium curves.}
\label{fig:totcollmnfld}
\end{figure}
\end{rem}
Next we describe the flow on the zero velocity manifold $\hat E(h, r,
\alpha)=0$, where 
the equations of motion reduce to
\begin{equation}\label{eq:motl=2hatE=0}
\dfrac{dr}{d\zeta}= 0, \quad
\dfrac{d\alpha}{d\zeta}= 0, \quad
\dfrac{d\psi}{d \zeta} = -\dfrac{r^2}{\sqrt \pi}\,
  \left[\dfrac{2\cos(\psi+\alpha)}{\sin 2\alpha} - \sin(\psi-\alpha)\right].
\end{equation}

\begin{lem}\label{thm:flussodegzerovelmnfld}
The flow on the total collision manifold is totally degenerate. More precisely
\begin{enumerate}
\item[(i)] $W^u(\vP_3) \equiv W^s(\vP_4)$; 
\item[(ii)] $W^u(\vP_4) \equiv W^s(\vP_3)$;
\end{enumerate}
where $\vP_3 \in \mathscr P_3$ and $\vP_4 \in \mathscr P_4$ are chosen
in such a way that the first two coordinates agree.
\end{lem}
\proof From Lemma \ref{thm:restpoints} we have
$\vP_4=\vP_3+(0,0,\pi)$. From the equations of motion
(\ref{eq:motl=2hatE=0}) follows that the flow lines on
$\hat{E}(h,r,\alpha)=0$ are straight lines parallel to the
$\psi$-axis. Therefore they are heteroclinic connections between each
point  $\vP_3 \in \mathscr P_3$ and the corresponding  $\vP_4 \in \mathscr P_4$.
\finedim 

By choosing $(\alpha, \psi)=(\pi/4, \pi+ \pi/4)$ or $(\alpha,
\psi)=(\pi/4, \pi/4)$ the last two equations in
\eqref{eq:mcgehee22trisdiel=2} are identically zero. This means that
the constant pair $(\alpha, \psi)=(\pi/4, \pi+ \pi/4)$ is a constant
solution of the subsystem obtained by projecting the system of ode's
onto the $(\alpha, \psi)$-plane. By summing up, the following result
holds.

\begin{figure}
\centering
{%
\includegraphics[width=0.440\textwidth]{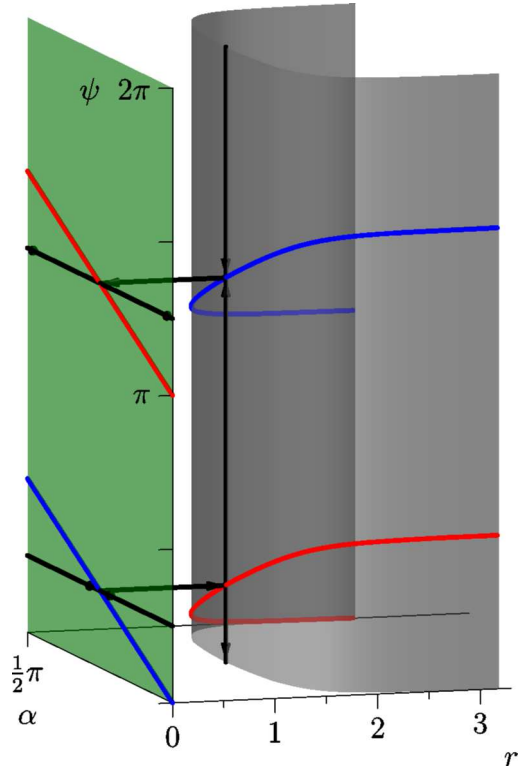}
}
\caption{Heteroclinic connections among the points $(0,\pi/4,\pi/4)$,
  $(0,\pi/4,\pi+\pi/4)$, $(\bar{r}, \pi/4, \pi/4)$, $(\bar{r}, \pi/4,
  \pi+\pi/4)$ where $\bar{r}$ is such that $\hat{E}(h,\bar{r},\pi/4)=0$.}
\label{fig:heteroclinic}
\end{figure}

\begin{prop}\label{thm:eterobetween}
  There exist two connecting orbits between the total collision and
  zero velocity manifold. More precisely
\begin{enumerate}
\item  The curve $\eta_{13}(\zeta)= \big(r(\zeta), \pi/4,  \pi/4\big)$ is an heteroclinic joining $\mathscr P_1$ and $\mathscr P_3$ for 
$r$ solution of 
\[
\dfrac{dr}{d\zeta}= \dfrac{1}{\sqrt \pi}\dfrac{r^3}{r^2+2} \, \hat E (h,r, \pi/4).
\]
\item the curve $\eta_{24}(\zeta)= \big(r(\zeta), \pi/4, \pi +
  \pi/4\big)$ is an heteroclinic joining $\mathscr P_4$ and $\mathscr
  P_2$ 
for $r$ solution of
\[
\dfrac{dr}{d\zeta}= -\dfrac{1}{\sqrt \pi}\dfrac{r^3}{r^2+2} \, \hat E (h,r, \pi/4). 
\]
\end{enumerate}
where $\hat E(h,r,\pi/4)= 2[h r^2 + 3(1-r^2\log r)]- 2r^2 \log 4$.
\end{prop}

\begin{rem}
  From Lemma \ref{thm:plcc} with $l=2$ it follows that the
  heteroclinic connections $\eta_{13}$ and $\eta_{24}$ are
  self--similar orbits where the vortices are located at the vertices
  of a square. Their projection on the shape sphere is a central
  configuration. Moreover, these are the only self--similar orbits:
  from \eqref{eq:mcgehee22trisdiel=2} it follows that no other solution curve
  has $\alpha(\zeta) = \pi/4$ when $\zeta$ ranges on an interval of non--zero length.
\end{rem}

\begin{rem}
  \label{rem:sundman}
  A classical result of Sundman \cite{sundman} proved for the
  three-body problem in three dimensions that an orbit ending in
  triple collision asymptotically approaches a central
  configuration. The validity of Sundman-type asymptotic estimates for
  collision solutions is established in a recent paper \cite[Theorem
  5, Example2]{gencoll} for a wide class of dynamical systems with
  singular forces, including the classical $N$-body problem,
  quasi-homogeneous and logarithmic potentials. Applying these results
  to our problem, we deduce that the only point of the total collision
  manifold which is (asymptotically) reachable by an orbit originating
  from an initial condition having $r>0$ is
  $(r,\alpha,\psi)=(0,\pi/4,5\pi/4)$, while
  $(r,\alpha,\psi)=(0,\pi/4,\pi/4)$, is the only point reachable extending the orbit
  backward in time. In this sense, the collision manifold, with the
  exception of these two points, is dynamically disconnected from the
  region $r>0$. 
\end{rem}
Figure \ref{fig:heteroclinic} shows graphycally the heteroclinic
connections between the invariant manifolds $r=0$ and $\hat{E}=0$.

\begin{rem}
  We stress that, while an initial condition on the heteroclinic
  $\eta_{13}$ is unable to complete the loop up to the total collision
  manifold (it just asymptotically reaches a fixed point on
  $\mathscr{P}_3$), the dynamics expressed in Cartesian coordinates by
  \eqref{eq:ham} does perform the loop in a finite time. In fact, the
  points on the curves $\mathscr{P}_3$ and  $\mathscr{P}_4$ of the zero velocity
  manifolds appear as rest points in McGehee coordinates just because
  the time transformation \eqref{eq:timescaling} is singular when
  $\hat{E}(h,r,\alpha)=0$. Their counterparts in Cartesian coordinates
  correspond to turning points where $\dot{\vq}(\sigma)=\zero$
\end{rem}

\section{Global dynamics in McGehee coordinates}

The typical orbit of equations \eqref{eq:mcgehee22trisdiel=2}
experiences a binary collision within a finite time, as we shall prove
in this section. The exceptions constitute a set of zero Lebesgue
measure. Moreover, there are no unbounded orbits, even when one
introduces generalized solutions that allow for an orbit to pass
through a binary collision.

\subsection{Non-existence of unbounded non-colliding trajectories}

Let us start by ruling out unbounded, collisionless orbits. In our
problem orbits escaping to infinity cannot be ruled out with a simple
argument based on the conservation of total energy $h$, because for
any given $h$ the distance $r$ from the center of vorticity may grow
without bounds while maintaining the same positive value of the
(rescaled) kinetic energy $\hat{E}$, as should be apparent from
Figures \ref{fig:energy} and \ref{fig:plotenrgiazerodiversil}.

\begin{defn}
  We say that a solution $\gamma(\zeta):=(r(\zeta), \alpha(\zeta),
  \psi(\zeta))$ of the system 
  \eqref{eq:mcgehee22trisdiel=2} is
  unbounded if
\[
\lim_{\zeta\to +\infty} r(\zeta) = + \infty
\]
\end{defn}

\begin{thm}\label{thm:nohyp}
  Equations \eqref{eq:mcgehee22trisdiel=2} do not allow for unbounded
  collisionless orbit.
\end{thm}
\proof Let us consider the system \eqref{eq:mcgehee22trisdiel=2};
eliminating the time variable $\zeta$ this system is equivalent to
\begin{equation}\label{eq:mcgeheesenzatempo}
 \left\{
\begin{array}{ll}
\dfrac{d\alpha}{dr} = -\dfrac{r^2+2}{r^3}\, \tan(\psi-\alpha)\\
\\
\dfrac{d\psi}{dr} = -\dfrac{r^2+2}{r \hat E}\left[\dfrac{2\cos(\psi+\alpha)}{\sin(2\alpha)\, \cos(\psi-\alpha)}- \tan(\psi-\alpha)\right]
\end{array}\right.
\end{equation}
By contradiction, we assume that there exists an unbounded motion $r
\mapsto (\alpha(r), \psi(r))$ starting from the initial condition
$\vh_0=(r_0, \alpha_0,\psi_0)$ and we suppose that the orthogonal
projection of $\vh_0$ onto the $(r, \alpha)$-plane is a point $\hat
\vh_0$ lying in the region between the straight lines $\alpha=0$,
$r=r_*$ and the zero set of the function $\hat E $, where $r_*$ is
such that $\hat E_0(r_*,\pi/4)=0$ (see Figure
\ref{fig:plotenrgiazerodiversil}). The case where $\hat \vh_0$ lies
between the straight lines $\alpha=\pi/2$, $r=r_*$ and the zero set of
the function $\hat E $ is analogous and will not be explicitly
addressed. In order to maintain $\hat{E}>0$, an unbounded motion must
satisfy
\begin{equation}\label{eq:unbounded_inequality}
 0 < \alpha(r)< g(r), \qquad \textrm{where}\quad  g(r):= \dfrac12\arcsin\left[\dfrac14 \, \exp\left(\dfrac{E_0(r)}{2\,r^2}\right)\right]
\end{equation}
for $E_0(r)=6(1-r^2\log r)$. 
We observe that 
\[
\lim_{r\to +\infty} E_0(r)/2r^2=\lim_{r\to +\infty}\log r^{-3}=-\infty,
\]
from which it follows that $g(r)\sim_{+\infty} 1/(8r^3)+ o(1/r^4)$. 
Thus $r \mapsto \alpha(r)$ decreases to zero at infinity while
keeping $r\mapsto g(r)$ as an upper bound. From the first 
equation in \eqref{eq:mcgeheesenzatempo} we have
\[
\dfrac{d\alpha}{dr}= - \dfrac{r^2+2}{r^2} \dfrac{\tan(\psi(r)-\alpha(r))}{r} \sim_{+\infty} - \dfrac{\tan(\psi(r)-\alpha(r))}{r}.
\]
and therefore 
\[
\lim_{r\to +\infty} \dfrac{\tan(\psi(r)-\alpha(r))}{r}=0,
\]
We conclude that the function $r \mapsto \alpha(r)$ should satisfy the limit equation 
\[
\dfrac{d\alpha}{dr}=- \dfrac{\tan(\psi(r))}{r}.
\]
From the second equation in \eqref{eq:mcgeheesenzatempo} we have
\[
\dfrac{d\psi}{dr} = -\dfrac{r^2+2}{r\hat E}\left[\dfrac{2\cos(\psi+\alpha)}{\sin(2\alpha)\cos(\psi-\alpha)} -\tan(\psi-\alpha)\right]\sim_{+\infty} 
-\dfrac{r}{\hat E}\left[\dfrac{1}{\alpha}-\tan \psi\right].
\]
From the upper bound  $r\mapsto \alpha(r) \in O(1/r^3)$ we also
have 
\[
\hat E(h,r,\alpha) \sim_{+\infty} - 6 r^2\log r + o(r^2\log r).
\]
Therefore an unbounded motion should satisfy the   
limit problem
\[
(\mathscr P_\infty)\ \quad   \left\{
\begin{array}{ll}
\dfrac{d\alpha}{dr}= -\dfrac{\tan\psi}{r}\\
\\
\dfrac{d\psi}{dr} = \dfrac{1}{6 r \log r}\left[\dfrac{1}{\alpha}-\tan\psi\right].                                                   
\end{array}\right.
\]
The first equation is equivalent to
\[
\psi(r)= -\arctan(r\alpha') 
\]
hence unbounded orbits should asymptotically satisfy the following
second order differential equation
\begin{equation}\label{eq:scalarnohyper}
- \dfrac{\alpha' +r \alpha''}{1+ r^2 \alpha'^2}= \dfrac{1}{6r\log r }\left[\dfrac{1+r\alpha\, \alpha'}{\alpha}\right].
\end{equation}
Multiplying the equation by $\alpha$ and taking into account that both
$\alpha$ and $\alpha'$ go to zero as $r\to +\infty$ it follows that
the above equation reduces to
\[
 \alpha \alpha''= -\dfrac{1}{6r^2 \log r}.
\]
Observing that $\alpha>0$, it follows that an unbounded solution
should be concave, which contradicts the inequality
\eqref{eq:unbounded_inequality}. \finedim

\begin{figure}[ht]
\centering
{%
\includegraphics[width=1.0\textwidth]{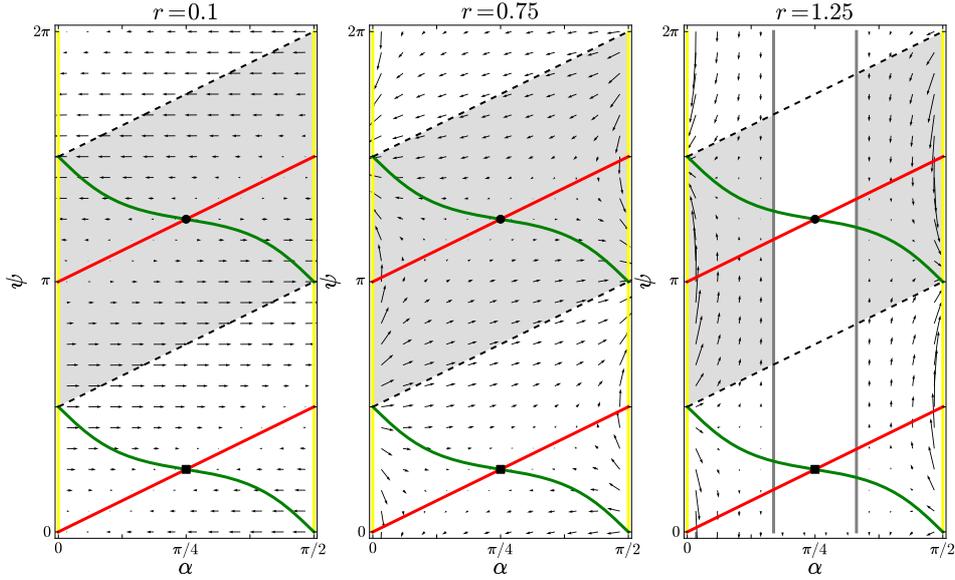}
}
\caption{
Projection of the vector field in equation \eqref{eq:mcgehee22trisdiel=2}
on the $(\alpha,\psi)$ plane at three different values of $r$. 
The black square and circle are, respectively the intersection with the
plane of the ejection and collision heterocline of Proposition
\ref{thm:eterobetween}.
The yellow vertical lines show the binary-collision planes at $\alpha=0$
and $\alpha=\pi/2$.
The red and green lines are, respectively, the places where 
$d\alpha/d\zeta=0$ and $d\psi/d\zeta=0$.
The vertical grey lines show the intersection of the zero velocity
manifold $\hat{E}=0$ with the plane $r=1.25$. The vector field is not
shown in the dynamically meaningless region between the two lines,
where $\hat{E}<0$. 
The regions shaded in light grey are those where $dr/d\zeta<0$.
Elsewhere it is $dr/d\zeta>0$, except along the black dashed lines,
where $dr/d\zeta=0$.
}
\label{fig:heteroclinicvectorfield2}
\end{figure}

\subsection{Generality of binary collisions}

\begin{thm}\label{thm:sibinary}
  Any orbit not fully contained in the total collision manifold or in
  the zero velocity manifold, or in the central manifold of the rest
  point $(r,\alpha,\psi)=(0,\pi/4,5\pi/4)$ or in the stable manifold
  of the rest points $\mathscr{P}_3$ experiences a binary
  collision. Moreover, the union of the orbits that do not experience
  a binary collision is a set of zero Lebesgue measure in $\mathbb{R}^3$.
\end{thm}
\proof The total collision manifold and the zero velocity manifold are
both two-dimensional, hence of zero Lebesgue measure. They are
invariant manifolds and the orbits lying on them do not experience
binary collisions, as illustrated in the previous section.

Of all the points on the total collision manifold, only
$(r,\alpha,\psi)=(0,\pi/4,5\pi/4)$ may be reached starting from $r>0$
(Remark \ref{rem:sundman}). The heterocline $\eta_{24}$ belongs to the
central manifold of $(r,\alpha,\psi)=(0,\pi/4,5\pi/4)$ associated to
the eigenvector $\vv_2$ of $\mathscr{P}_2$ (Lemma
\ref{thm:spectrum1}). Owing to the non-uniqueness of the central
manifold, other orbits are possible that become tangent to the
heterocline $\eta_{24}$ as $r\to 0$.  Defining $\phi=\psi-\alpha-\pi$,
$\beta = \alpha -\pi/4$, and approximating equations
\eqref{eq:mcgehee22trisdiel=2} at first order in $\beta$, $\phi$, the
dynamics of orbits arbitrarily close to the heterocline is described by
\begin{equation}\label{eq:mcgeheeD2_linearization}
\left\{
\begin{array}{ll}
\dfrac{dr}{d\zeta} = - \hat{E}\,\dfrac{r^3}{r^2+2} \\
\\
\dfrac{d\beta}{d\zeta} = - \hat{E}\,\phi\\
\\
\dfrac{d\phi}{d \zeta} = -r^2 (3\phi+4\beta).
\end{array}\right.
\end{equation}
where we have omitted the irrelevant $\pi^{-1/2}$ factors. Here
$(\beta, \phi)=(0,0)$ is the heterocline, and $r$ is a monotonically
decreasing function of $\zeta$.  At each fixed value $r>0$ this system
of equations has a contracting and an expanding direction on the
$(\beta, \phi)$-plane, associated, respectively, to the eigenvalues
$\lambda_1=-r(\sqrt{16\hat{E}+9r^2}+3r)/2$ and
$\lambda_2=r(\sqrt{16\hat{E}+9r^2}-3r)/2$. The contracting direction
is given by the vector
$(r,\beta,\phi)=(0,2\hat{E},r(\sqrt{16\hat{E}+9r^2}+3r))$, therefore
the center manifold is a two-dimensional surface that becomes tangent
to $(r,\alpha,\psi)=(r,\alpha,\alpha+\pi)$ at $(\alpha,\psi)=(\pi/4,5\pi/4)$ as $r\to 0$.

Another set of orbits that do not experience binary collisions are the
stable manifolds of the points of $\mathscr{P}_3$ associated to the stable eigenvector 
$\vu_2$ (Lemma \ref{thm:spectrum2}). These points are normally
hyperbolic, and the union of their stable manifolds is, again, a
smooth two dimensional surface. See \cite{shub} for details.

Finally, from equations \eqref{eq:mcgehee22trisdiel=2}, in the
rectangle $\vR=(0,\pi/2)\times[0,2\pi]$ of the $(\alpha, \psi)$-plane,
we have the following properties of the orbits (see also Figure
\ref{fig:heteroclinicvectorfield2}):
\begin{itemize}
\item $\alpha$ increases (with respect to time $\zeta$) in the region
  of $\vR$ where
\[
\alpha < \psi <\alpha+ \pi
\] 
and decreases otherwise;
\item $\psi$ increases (with respect to time $\zeta$) in the region of $\vR$ where
\[ 
\arctan(m(\alpha))<\psi<\arctan(m(\alpha))+\pi.
\]
and decreases otherwise. 
\item $r$ decreases (with respect to time $\zeta$) in the region of $\vR$ where
\[ 
\alpha+\pi/2 < \psi <\alpha+ 3\pi/2
\]
and decreases otherwise. 
\end{itemize}
The qualitative behavior of the vector field then rules out limit
cycles or other invariant structures.  Theorem \ref{thm:nohyp} shows
that collisionless, unbounded orbits are impossible.  It then follows
that orbits not on the above mentioned invariant manifolds, whose
union is a set of zero Lebesgue measure, must reach a binary
collision.\finedim


\subsection{Generalized solutions: a boundedness result}

As a consequence of the fact that almost every orbit experiences a
binary collision, we are lead to ask if it is possible to continue the
solution after a binary collision in some {\em natural\/} manner. This
is classical question in the context of $n$-body gravitational
problems and it goes under the name of {\em regularization of
  collisions.\/} Motivated by the recent paper \cite{cate} we shall
introduce the notion of {\em generalized solutions\/} which, very
roughly, are the continuation of the classical solutions after the
binary collision in a suitable manner. In this section all
``collisions'' should be intended as ``binary collisions'', thus excluding
the case of total collapse.
Before proceeding further, we firstly recall some facts proven in
\cite{cate} in the case of weak central forces.

\subsubsection*{One-center problem with logarithmic potential}

Let us consider the dynamical system associated with the conservative
central weak force arising by a logarithmic potential having the 
singularity at the origin; i.e. the Cauchy problem:
\begin{equation}\label{eq:centralforce}
 \left\{
\begin{array}{ll}
\ddot \vq = \partial_\vq U(\vq)\\
\\
(\vq(0), \dot \vq(0))=(\vq_0, \vp_0) \in (\R^2\setminus\{\zero\} \times \R^2)
\end{array}\right.
\end{equation}
where $U(\vq):=-\log(\|\vq\|)$. A {\em classical solution\/} does not
cross the singularity of the force, i.e. is a path $\vq \in \mathscr
C^2(T, \R^2\setminus\{\zero\})$ where $T$ denotes the maximal interval
of existence. For the classical $n$-body problem, De Giorgi in
\cite{Degiorgi96} proposed in to consider the smoothing of the
potential as a regularization technique. An exhaustive analysis in
the case of homogeneous potentials, including the classical Keplerian
potential, was performed by Bellettini, Fusco e Gronchi in
\cite{bellettini}. Following the idea of De Giorgi, the authors in
\cite{cate} removed the singularity at $\vq=0$ by smoothing the
potential as follows:
\[
 U(\vq; \varepsilon)= \log \sqrt{\|\vq\|^2 + \varepsilon^2}
\]
and by considering the regularized problem
\begin{equation}\label{eq:centralforceepsilon}
 \left\{
\begin{array}{ll}
\ddot \vq = \partial_\vq U(\vq;\varepsilon)\\
\\
(\vq(0), \dot \vq(0))=(\vq_0, \vp_0) \in (\R^2 \times \R^2)
\end{array}\right.
\end{equation}
Unlike in \eqref{eq:centralforce} the differential equation in
\eqref{eq:centralforceepsilon} is no longer singular, so the initial
value problem admits a global smooth solution $\vq \in \mathscr
C^\infty(\R, \R^2)$ for every choice of the initial value $(\vq_0,
\vp_0)$. Let $B_0(\bar R)$ be a ball of radius $\bar R$ centered at
the origin of the configuration space, and let $\mathscr S\subset \R^2
\times \R^2$ be the set of initial conditions within the ball leading
to collision for the problem with $\varepsilon=0$. For every $\bar \nu
\in \mathscr S$ let $\vq_{\bar \nu} \in \mathscr C^2(T, \R^2)$ be the
collision solution where $T$ denotes the maximal interval of
existence. Denoting with $\vq_{\varepsilon, \nu}$ the solution with
initial data $\nu$ and studying the asymptotic limit 
$(\varepsilon , \nu)\to (0, \bar \nu)$, the authors in \cite{cate}
introduced the following two notions of regularization.
\begin{defn}\label{def:weakreg}
  We say that the problem \eqref{eq:centralforce} is {\em weakly
    regularizable\/} via smoothing of the potential in the ball
  $B_0(\bar R)$ if, for every $\bar \nu \in \mathscr S$ there exist
  two sequences $(\varepsilon_k)_k, (\nu_k)_k$ tending to $0$ and
  $\bar \nu$ respectively such that there exists
\[
 \lim_{k\to \infty} \vq_{\varepsilon_k, \nu_k} = \bar \vq
\]
and the flow 
\[
 \tilde \vq_\nu (t)=\left\{
\begin{array}{ll}
\vq_\nu(t) & \nu \notin \mathscr S\\
\bar \vq (t) &\nu \in \mathscr S
\end{array}\right.
\]
is continuous with respect to $\nu$ (i.e. the initial point).
\end{defn}
\begin{defn}\label{def:strongreg}
 The singular one center problem \eqref{eq:centralforceepsilon} is
 strongly regularizable via smoothing of the potential if there 
exists $\bar R$ such that for every $\bar \nu \in \mathscr S$ the limit 
\[
\lim_{(\varepsilon, \nu)\to (0, \bar \nu)} \vq_{\varepsilon, \nu} = \bar \vq 
\]
exists and the flow 
\[
\tilde \vq_\nu(t)=\left  \{
\begin{array}{ll}
\vq_\nu(t) & \nu \notin \mathscr S\\
\bar \vq (t) &\nu \in \mathscr S
\end{array}\right.
\]
is continuous with respect to $\nu$.
\end{defn}
The authors in \cite{cate} then prove that
\begin{prop}
  The logarithmic one center problem is globally regularizable via
  smoothing the potentials according to Definition
  \ref{def:strongreg}.
\end{prop}
As a consequence of this result, 
it follows that it is possible to continue the solutions after the
collision by transmission, in the following sense.
\begin{defn}\label{def:transsolutions}
  Let $t \mapsto \vq_0(t)$, $t \in [0, T_0)$ be a collision path, and
  $T_0$ the collision instant. Define the {\em transmission
    solution\/} $\bar \vq_0$ for $t \in [0, 2T_0]$ as follows:
\[
\left\{
\begin{array}{ll}
 \bar \vq(t)= \vq_0(t) & t \in [0, T_0]\\
\\
\bar \vq(t)=- \vq_0(2T_0-t) & t \in[T_0, 2T_0]
\end{array}\right. 
\]
\end{defn}
This is the right definition to set in order to have the continuity
with respect to the initial data of the flow obtained by replacing
the collision solution $\vq_0$ with the transmission solution.

\subsubsection*{Local regularization of collisions and generalized solutions}

We apply the above results to the solutions of our problem that
experience binary collisions. The presence of other vortices, of
course, complicates the picture, but it should be intuitively clear
that any binary collision, locally, is just a central problem, as we
show formally below.


In Cartesian coordinates, let us consider instead of the singular
potential defined in \eqref{eq:potenziale}, the non-singular one given
by
\begin{equation}\label{eq:potenziale-cate}
U(\vq_0; \varepsilon) = -\sum_{g\in D_l \smallsetminus \{1\}}
\log\left(\sqrt{\left\| \vq_0 - g \vq_0 \right\|^2 +\varepsilon^2}\right),
\end{equation}
that substitutes (\ref{eq:potenziale}) in the equations of motion
(\ref{eq:ham}).
For the Klein group $D_2$ (see section \ref{sec:stabunstab}), the above
expression reduces to
\begin{equation}\label{eq:Klein-potential-cate}
U(q_1,q_2;\varepsilon) = \frac{1}{2} \log\left(
(4q_1^2+\varepsilon^2)(4q_2^2+\varepsilon^2)(4\left\|q_0\right\|^2+\varepsilon^2)\right).
\end{equation}
We recall that a simultaneously binary collision in our problem
corresponds to a solution $\bar\vq$ of the Newton's equations
\eqref{eq:Newton}
\[
\Gamma \ddot \vq = \dfrac{\partial U}{\partial \vq}.
\]
for $U(\vq):= -\log(8 \|\vq\|q_1 \, q_2)$ such that at some time
instant $\sigma_*$ one and only one coordinate of the point $\bar
\vq(\sigma_*):=(q_1(\sigma_*), q_2(\sigma_*))$ is zero. This means
that the support of the solution intersects one of the coordinate axis
$q_1$ or $q_2$.  Without loss of generality we can assume that at
some instant $\sigma_*$ we have $q_2(\sigma_*)=0$. Thus there exists
$\delta:= \delta(\sigma_*)>0 $ such that
\[
-\log(q_2\, K_1) \leq U(\vq) \leq -\log(q_2\, K_2),
\]
where $K_1:= 8 \, (C+ \delta) ^2$ and $K_2:= 8\, (C-\delta)^2$ and
$C:=q_1(\sigma_*)$. Therefore, locally in the neighborhood of a
singularity, the potential is, up to a constant, the central
logarithmic potential described in the paragraph above.
\begin{rem}
  By continuing the solution after a binary collision as described
  above we can replace the colliding trajectory, which does not exists
  beyond $\sigma_*$, with a transmission solution that exists also at
  later times. However at some instant this new trajectory could
  experience another binary collision that can, again, be extended by
  transmission, and so on for any number of collisions. Although we
  could construct arbitrarily long solutions containing an infinite
  number of collisions, the above arguments do not guarantee the
  uniform convergence of the generalized solutions. The extension of a
  generalized solution up to an infinite time (if possible) would
  require  a much deeper and careful analysis than that presented above.
\end{rem}
\begin{defn}\label{def:gencoll1}
  A {\em generalized solution\/} $\vq \in \mathscr C^2(\hat J) \cap
  \mathscr C^0(J)$ of the problem is 
  a classical solution
  (``classical'' meaning a $\mathscr C^2$ curve that satisfies the
  equations\eqref{eq:Newton} pointwise) on the time interval
  $\hat
  J=J\setminus Y$, where $Y:=\{\sigma_1, \dots. \sigma_n\}$ is a finite
  set of times. At those times the solution $\vq$ experiences a binary
  collision and is continued by transmission as explained above.
\end{defn}

Before proceeding further we translate the notion of generalized
solution to the McGehee setting. In these coordinates, a solution
experiences a binary collision when it crosses either the plane
$\alpha=0$ or the plane $\alpha=\pi/2$, corresponding, respectively,
to collisions with $q_2(\sigma_*)=0$ and  $q_1(\sigma_*)=0$ in
Cartesian coordinates. 
The graph of $r\mapsto \alpha(r)$ on the $(r,\alpha)$ plane for a
transmission orbit colliding at $r=r_*$ is locally symmetric with respect
to $r_*$. The graph of $r\mapsto \psi(r)$ has a jump discontinuity of
$\pi$ at $r=r_*$.


\subsubsection*{Boundedness of the generalized orbits}

The equations of motion given in \eqref{eq:mcgehee22trisdiel=2}
may be written as follows:
\begin{equation}\label{eq:mcgehee22trisdiel=2new}
\left\{
\begin{array}{ll}
 \dfrac{d\alpha}{dr}= \dfrac{r^2+2}{r^3}\,\tan(\psi-\alpha)\\
\\
\dfrac{d\psi}{dr}= -\dfrac{r^2+2}{r\, \hat E}\left[\dfrac{2\cos(\psi+\alpha)}{\sin (2\alpha)\, \cos(\psi-\alpha)}
-\tan(\psi-\alpha)\right]\\
\\
\dfrac{d\psi}{d \alpha} = \dfrac{- r^2}{\hat E}\left[\dfrac{2\cos(\psi+\alpha)}{\sin (2\alpha)\, \sin(\psi-\alpha)}-1\right]
\end{array}\right.
\end{equation}
As observed in Theorem \ref{thm:nohyp}, if $r\to +\infty$ then
necessarily $\alpha\to0$, otherwise the constraint $\hat E >0$ would be
violated. In the following we shall distinguish between two cases:
\begin{enumerate}
\item if $\alpha(r) \in O(1/r^3)$ for $r\to \infty$  then $\hat E_0(r, \alpha)\approx c$ for some positive constant $c$;
\item if $\alpha(r) \in o(1/r^3)$ for $r\to\infty$ then $\hat E_0(r, \alpha)\approx -2r^2\log \alpha$.
\end{enumerate}
{\em First case.\/} For fixed $\psi$ we have that
 \[\tan(\psi-\alpha)=  \dfrac{\tan\psi }{1+\alpha\tan\psi} + o(\alpha),\]
\[
\cos(\psi-\alpha)= \cos\psi+ \alpha \sin\psi + o(\alpha),\qquad 
\cos(\psi +\alpha)= \cos\psi - \alpha \sin\psi+o(\alpha),
\] 
and hence the system given in \eqref{eq:mcgehee22trisdiel=2new} reduces to 
 \begin{equation}\label{eq:asint1}
\left\{
\begin{array}{ll}
\dfrac{d\alpha}{dr}= \dfrac{1}{r}\,\tan(\psi)\\
\\
\dfrac{d\psi}{dr}= -r \left[\dfrac{1}{\alpha}-2\tan\psi\right]\\
\\
\dfrac{d\psi}{d \alpha} = -r^2\left[\dfrac{1-2\alpha \tan\psi}{\alpha\tan\psi}\right].
\end{array}\right.
\end{equation}
\begin{lem}\label{thm:lemmamonotone1}
  Let us consider an orbit starting at the point $P_0=(r_0, \alpha_0,
  \psi_0)$ and let us 
  also assume that $\tan \psi$ is bounded. Then
  the projection of the orbit onto the $(\alpha,\psi)$-plane
  intersects the line $\alpha=0$ for $r>r_0$ if $\tan \psi<0$ and for
  $r <r_0$ otherwise.
\end{lem}
\proof The result immediately follows from the first equation in
\eqref{eq:asint1}.\finedim\\
In order to show that generalized solutions are bounded, we must
consider the following three cases: 
\begin{enumerate}
 \item[(i)]
If $\alpha \tan\psi \in O(1)$  then the system \eqref{eq:asint1} reduces to
\begin{equation}\label{eq:asint2}
\dfrac{d\alpha}{dr}\sim \dfrac{1}{\alpha\,r},\qquad
\dfrac{d\psi}{dr}\sim \dfrac{r}{\alpha}, \qquad
\dfrac{d\psi}{d \alpha} \sim r^2
\end{equation}
\item[(ii)]If $\alpha \tan\psi \in o(1)$ then we have 
\begin{equation}\label{eq:asint3}
\dfrac{d\alpha}{dr}\sim \dfrac{\tan\psi}{r},\qquad
\dfrac{d\psi}{dr}\sim \dfrac{r}{\alpha}, \qquad
\dfrac{d\psi}{d \alpha} \sim \dfrac{r^2}{\alpha \tan \psi}.
\end{equation}
\item[(iii)]If $\alpha \tan\psi$ is unbounded then 
\begin{equation}\label{eq:asint4}
\dfrac{d\alpha}{dr}\sim \dfrac{\tan\psi}{r},\qquad
\dfrac{d\psi}{dr}\sim r \tan \psi, \qquad
\dfrac{d\psi}{d \alpha} \sim r^2.
\end{equation}
\end{enumerate}
An interesting consequence of Lemma
\ref{thm:lemmamonotone1} is that a
generalized solution may temporarily exhibit a 
monotonicity (with respect to $r$) of the collision
points. More precisely, taking into account the above estimates for 
$d\psi/dr$, the following result holds.
\begin{lem}\label{eq:lemmamonotone2}
  Under the assumptions of lemma \eqref{thm:lemmamonotone1} it follows
  that the projection of the orbit onto the $(\alpha,\psi)$-plane
  intersects the line $\alpha=0$ at the instants $(r_i)_{i \in \N}$
  with $r_i <r_{i+1}$ if $\tan \psi<0$ and $r_i >r_{i+1}$ otherwise.
\end{lem}
However, the same asymptotic estimates above imply that a
monotonic sequence of collisions cannot be arbitrarily long.
Otherwise stated, $\tan \psi$ is not bounded along the whole orbit.
\begin{lem}
  If $\alpha(r) \in O(1/r^3)$ unbounded generalized solutions do not exist.
\end{lem}
\proof By contradiction, we observe that in an unbounded solution
$\gamma$ the first component $r\mapsto \alpha(r)$ tends either to $0$
or to $ \pi/2$ for $r \to \infty$.  Without loss of generality we
discuss only the case $\alpha=0$.  We also observe that along an
unbounded generalized solution $\frac{d\alpha}{dr}\to 0$ when $r\to
+\infty$.  We are also 
assuming that $\alpha(r) \in O(1/r^3)$. This
immediately rules out case (i) above, as it would imply
$\frac{d\alpha}{dr}\to + \infty$.  In case (ii) $\psi$ is unbounded
because of $\alpha(r) \in O(1/r^3)$.  For the same reason, in case
(iii) $\tan(\psi)$ is unbounded.  Therefore there exists a $r_*$ such
that $\psi(r_*) =\pi/2 + k\pi$ for some $k \in \Z$.  For $r>r_{*}$ the
function $r\mapsto \alpha(r)$ decreases and the monotonically growing 
sequence of collisions stops. \finedim

{\em Second case.\/} In this case the system
\eqref{eq:mcgehee22trisdiel=2new} reduces to
 \begin{equation}\label{eq:asint12}
\left\{
\begin{array}{ll}
\dfrac{d\alpha}{dr}= \dfrac{1}{r}\,\tan(\psi)\\
\\
\dfrac{d\psi}{dr}= \dfrac{1}{2r\log\alpha} \left[\dfrac{1}{\alpha}-2\tan\psi\right]\\
\\
\dfrac{d\psi}{d \alpha} = -\dfrac{1}{\log\alpha}\left[\dfrac{1-2\alpha \tan\psi}{\alpha\tan\psi}\right].
\end{array}\right.
\end{equation}
By arguing exactly as in the previous case, in this case is also
possible to show that the following result holds.
\begin{lem}
  If $\alpha(r) \in o(1/r^3)$ unbounded generalized solutions do not exist.
\end{lem}

\section{Global dynamics in physical coordinates}

From now own, when we refer to {\em solutions} we mean
{\em generalized solutions} that may go through binary collisions.

The transformations linking the physical coordinates $(q_1, q_2)$ and
momenta $(p_1, p_2)$ to the McGehee coordinates $(r, \alpha, \psi)$,
for a given value of the total energy $h$, are
\begin{equation}\label{eq:trasftornaindietropos}
\left\{
\begin{array}{ll}
q_1=r\, e^{-1/r^2}\cos\alpha\\
q_2=r\, e^{-1/r^2} \sin\alpha
\end{array}\right.
\end{equation}
\begin{equation}\label{eq:trasftornaindietromom}
\left\{
\begin{array}{ll}
p_1=z_1/r=\sqrt{\hat E\,\pi}/r\,\cos\psi\\
p_2=z_2/r=\sqrt{\hat E\,\pi}/r\,\sin\psi
\end{array}\right.
\end{equation}
We recall that the various time scalings introduced along calculations
have the effect that both the invariant manifold corresponding to the
total collapse and that corresponding to zero velocity are reached
asymptotically as the rescaled time goes to infinity. 

We further observe that not all the rest points on the total collision
manifold have a physical meaning. Only the points $(r,
\alpha,\psi)=(0, \pi/4, \pi/4)$ and $(r, \alpha,\psi)=(0, \pi/4,
5\pi/4)$ are physically relevant. They correspond to central
configurations in which each vortex lies on the vertex of a square,
having the barycenter at the origin and the vertices on the
bisectrices of the coordinate axis. Thus the two heteroclinic
connections of Proposition \ref{thm:eterobetween} joining the total
collision manifold and the zero velocity manifold correspond,
respectively to the homographic ejection orbit from total collapse to
the zero velocity manifolds and to the homographic collision from the
zero velocity to the total collapse manifolds.

\begin{figure}
\centering{%
\includegraphics[width=\textwidth]{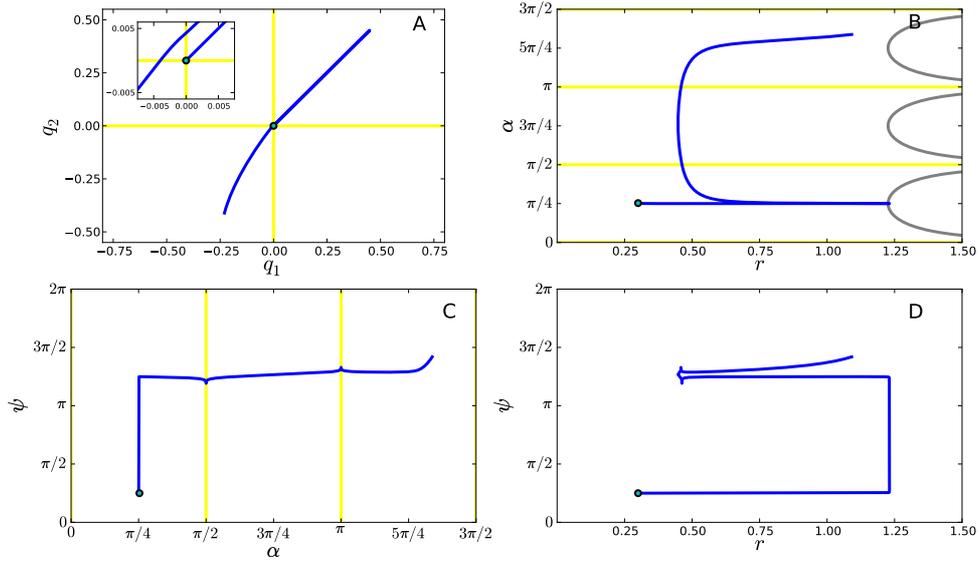}
}
\caption{A near-ejection orbit starting at $(r, \alpha,\psi)=(0.3,
  \pi/4+10^{-2}, \pi/4)$ with $h=0$. The initial conditions are
  transformed to physical coordinates using equations
  (\ref{eq:trasftornaindietropos}) and
  (\ref{eq:trasftornaindietromom}), then a numerical integration of
  the equations (\ref{eq:ham}) with the potential
  (\ref{eq:Klein-potential-cate}) and $\epsilon=10^{-6}$ generates the
  orbit.  Panel A): the orbit in physical coordinates, projected on
  the $(q_1,q_2)$ plane. The inset magnifies the region close to the
  origin (corresponding to the total collapse). Panels B), C), D): the
  orbit transformed back to McGehee coordinates, projected
  respectively on the $(r,\alpha)$, $(\alpha, \psi)$, and $(r, \psi)$
  planes. For clarity, after each binary collision the orbit is extended with
  continuity into the adjacent domain, rather than into the
  fundamental one, which would introduce visually bothersome jumps in the
  $\psi$ coordinate.  The circle shows the position
  of the initial condition. The yellow lines at $\alpha = 0, \pi,
  \cdots$ correspond to binary collisions at $q_1=0$, those at $\alpha
  = \pi/2, 3\pi/2, \cdots$ to binary collisions at $q_2=0$. The grey
  lines in panel B) are the points with $\hat{E}=0$.  }
\label{fig:ejection}
\end{figure}

This is illustrated in Figure \ref{fig:ejection} by a numerical
solution of the equations of motions in physical coordinates using the
non-singular potential (\ref{eq:Klein-potential-cate}) for an initial
condition close to the ejection heterocline. Initially the solution
follows closely the heteroclinic cycle shown in figure
\ref{fig:heteroclinic}. Eventually, the solution leaves the collision
heterocline before reaching the total collapse (which is a saddle
point) and it is subject to a sequence of two binary collisions close
to total collapse.

\begin{figure}
\centering{%
\includegraphics[width=\textwidth]{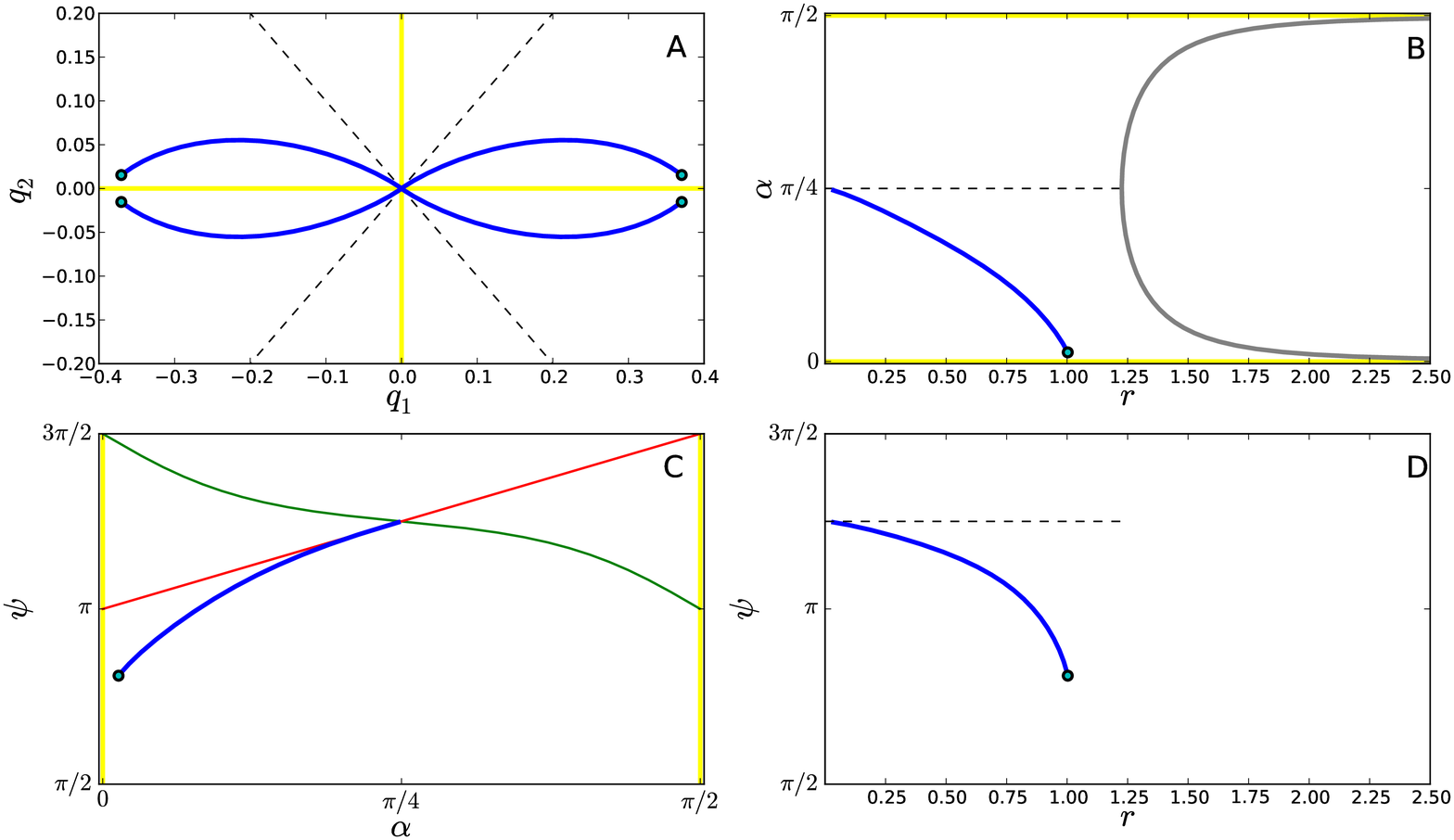}
}
\caption{An orbit as close as numerically possible to one lying on the
  center manifold of the total collision point
  $(r,\alpha,\psi)=(0,\pi/4,5\pi/4)$. The orbit starts approximately
  at $(r, \alpha,\psi)=(1.00264534\cdots,0.04149204\cdots,
  2.54344846\cdots])$ with $h=0$. The dashed lines are the homothetic
  orbit corresponding to the heterocline $\eta_{24}$. The red and
  the green curve in panel C) are, respectively, the projection on the
  $(\alpha,\psi)$ plane of the lines of rest points $\mathscr{P}_2$
  and $\mathscr{P}_4$. Details of the numerical solution
  and panel description as in figure \ref{fig:ejection}.
  }
\label{fig:pseudolemniscata}
\end{figure}

In Figure \ref{fig:pseudolemniscata} we illustrate a case
complementary to the above one. We show a numerical approximation to
an orbit lying in the center manifold of the total collision point at
$(r,\alpha,\psi)=(0,\pi/4,5\pi/4)$. The orbit approaches the
homothetic configuration (which, in McGehee coordinates is the
heteroclinic connection $\eta_{24}$) with a projection on the
$(\alpha,\psi)$ plane which is tangential to the line of rest points
$\mathscr{P}_2$. In Cartesian coordinates our choice of initial
conditions corresponds to the four vortices initially close in pairs,
as if they were just emerging from a binary collision, but with the
special choice of the momenta that sends them on a trajectory tangent
to the homothetic solution.


Summarizing we have
\begin{thm}\label{thm:2A}
  The only homothetic ejection/collision orbit from/to total collapse
  is that homothetic to the planar central configuration, which is bounded.
\end{thm}
\proof The existence and boundedness of the ejection/collision orbit
homothetic to the planar central configuration is a direct consequence
of Proposition \ref{thm:eterobetween}. For a hypothetical solution
where the position of the vortices is at all times a homothetic
transformation of the vertices of a given rectangle, then the momentum
of each vortex would be parallel to the line joining the vortex and
the origin.  However, it is clear that only at $\alpha=\pi/4$ all the
forces are balanced, therefore we exclude other solutions.  \finedim
\begin{rem}
  Recalling from Sundman type estimates proven in \cite{gencoll} that, of all the
  points on the collision manifold, those corresponding to a
  homothetic ejection/collision (namely $(\alpha, \psi)=(\pi/4,\pi/4)$
  and $(\alpha, \psi)=(\pi/4,5\pi/4)$) are dynamically reachable from
  outside the manifold, then the other rest points on the collision
  manifold and the heteroclinic connections between them are
  dynamically meaningless.
\end{rem}

\begin{thm}\label{thm:3A}
  Each solution different from the homothetic central configuration
  experiences a binary collision in finite time. Moreover, any
  solution starting arbitrarily close to the homothetic ejection
  solution has an orbit that may not reach the upper bound of the
  homothetic orbit before experiencing a binary collision.
\end{thm}
\proof This is a direct consequence of the theorem \ref{thm:sibinary}.
\finedim

\begin{figure}
\centering{%
\includegraphics[width=\textwidth]{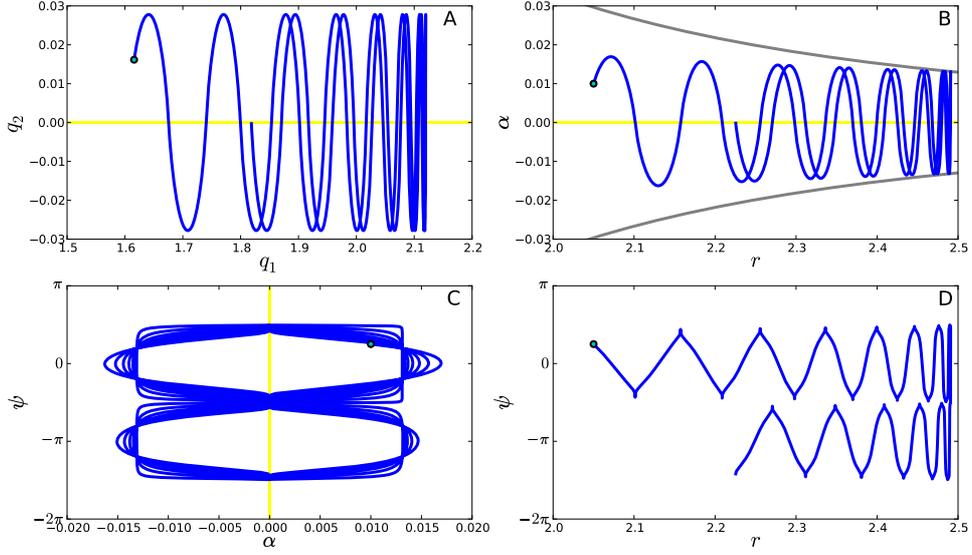}
}
\caption{A sequence of binary collisions moving along the positive
  $q_1$ semi-axis with a turning point. The orbit starts at $(r,
  \alpha,\psi)=(2.05, 0.01, \pi/4)$ with $h=0$. Details of the
  numerical solution and panel description as in figure
  \ref{fig:ejection}. }
\label{fig:turning-point}
\end{figure}
\begin{figure}
\centering{%
\includegraphics[width=\textwidth]{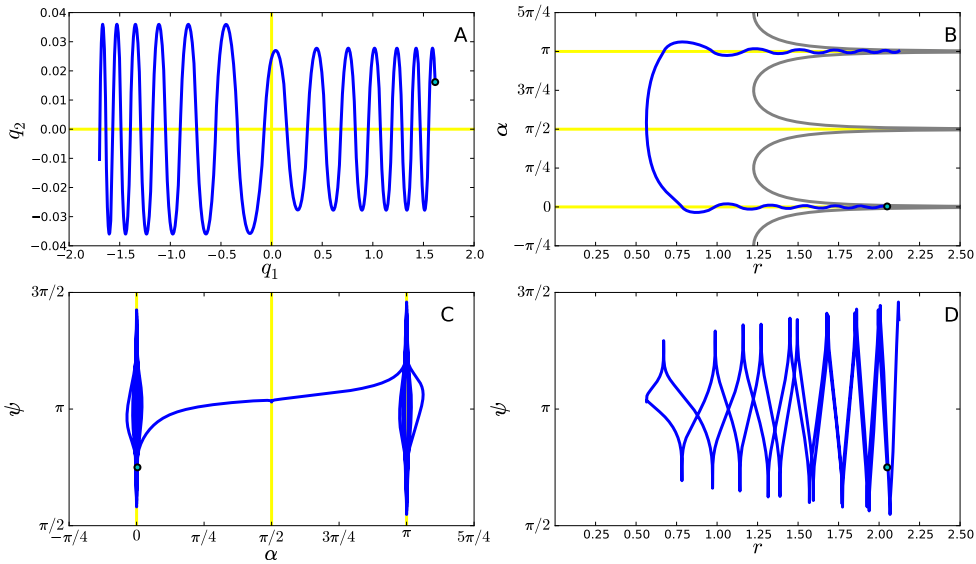}
}
\caption{A sequence of binary collisions moving from the positive to
  the negative $q_1$ semi-axis. The orbit starts at $(r, \alpha,\psi)=(2.05,
   0.01, 3\pi/4)$ with $h=0$. Details of the numerical solution and
   panel description as in  figure \ref{fig:ejection}. }
\label{fig:transmission}
\end{figure}

When the four vortices are arranged at the vertex of a rectangle with
a large ratio of side lengths, then there are solutions in which the
pairs of vortices undergo repeated binary collisions along the same
axis. In particular, for some of these solutions the sequence of
collisions initially moves away from the origin of the Cartesian
axes. However, unbounded solutions are not possible, and eventually
the sequence reaches a turning point and returns towards the
origin. An example of this dynamical behavior is shown in figure
\ref{fig:turning-point}.  When seen in the McGehee coordinates, these
are solutions advancing (up to a turning point) into the narrow funnel
between the zero velocity manifolds of adjacent domains. In figure
\ref{fig:transmission} we show a numerical solution where a sequence
of binary collisions approaching the origin from the right undergoes a
single binary collision along the $q_2$ axis and the continues moving
leftward on the negative $q_1$ semi-axis, although with a markedly
different amplitude. Eventually, this it will reach a turning point,
and return towards the origin. Solutions that oscillate along an axis
bouncing back and forth between turning points of opposite sign are
not uncommon. However, we were unable to maintain the oscillations,
which are aperiodic, for arbitrarily long times. Eventually the orbit
spends some time winding closely around the origin at distances
corresponding to values of $r$ such that $\hat{E}(h,r,\alpha)>0$ for
any $\alpha$. After these complicated and chaotic-looking transients
the orbit may resume its oscillations along one (not necessarily the
same) axis.

We may formalize these findings as follows.
\begin{thm}\label{thm:4A}
  If the initial configuration is a rectangle with large enough aspect
  ratio, the configuration evolves by passing through an arbitrary number of
  simultaneous binary collisions between the pairs of closest vortices.
\end{thm}
\begin{thm}\label{thm:5A}
The projection of any solution in the configuration space is bounded. 
\end{thm}
\proof This is a direct consequence of theorem
\ref{thm:nohyp}.\finedim

\begin{thm}\label{thm:6A}
  The set of initial conditions leading to total collapse or
  leading in an infinite amount of time to the outermost boundary region
  where the motion is allowed has zero Lebesgue measure.
\end{thm}
From the physical point of view the outermost boundary region among
the others depends on the total energy $h$ of the system and it is
given by the hypersurface having the property that any solution of the
Newton equations after reaching with zero velocity, the solution orbit
falls down.

Finally, we may go back to the initial interpretation of the equations
of motion (\ref{eq:Newton}) as steady solutions of the equations
(\ref{eq:scroedingernonlineare}) for infinitely tall, nearly parallel,
vortex filaments in a three-dimensional space. In that setting our
time-like variable really is the arc-length parameter along the vortex
filament. From this point of view an orbit of our dynamical system
describes the mutual position of four steady vortex filaments at
various heights. The points of intersection of the
vortex filaments with an horizontal plane, and the inclination of
tangent to the filaments at the intersection points are the
parameters that completely determine the shape of the whole filament.
In particular, some of our findings about the
equivariant solutions subject to the $D_2$ symmetry may be recast as
follows:
\begin{enumerate}
\item If the vortex filaments join together into a total collapse, or
  separate from collapse, then they do so by approaching a square
  configuration;
\item the distance among the vortex filaments is bounded at all
  heights;
\item binary collisions of the vortex are generic, that is, the set of
  intersection points and tangents generating non-colliding filaments
  has zero Lebesgue measure.
\end{enumerate}

\section{Further perspectives and closing remarks}

Among all the solutions of a system of vortex filaments some special
and very important class is represented by the so-called {\em
  helicoidal solutions\/}, namely solutions of
\eqref{eq:scroedingernonlineare} of the form.
\[
\Psi(\sigma,t)=e^{i\nu t}\vq(\sigma)
\]
for some complex valued function $\vq:\R\to\C$. In the special case in
which $\vq(\sigma)=\lambda(\sigma)\vxi$ it is called a {\em
  homographic solution\/} and if if $\lambda(\sigma)=e^{i\nu\sigma}$
we shall refer as a {\em relative equilibrium.\/}
We observe that, for fixed frequencies $\omega, \nu$ the function
\[
\Psi(\sigma, t)= e^{i\omega \sigma}\, e^{i \nu t}\vxi
\]
only depends on the choice of the configuration $\vxi$. Moreover every
filament has the same shape of helicoidal type and the time evolution
corresponds to a rotation equal for each vortex filament of the
initial configuration. Therefore the solution $\Psi$ is periodic in
$\sigma$ and constant in time.
\begin{enumerate}
\item An interesting question that should be addressed is to consider 
helicoidal solutions instead of stationary solution. In this case the 
dihedral potential we need to add a quadratic term.  However we think that 
our techniques should work also in this case. 
\item Of some interest in the applications is the study of a singular logarithmic type potential on a 
more general surface without boundary. 
\end{enumerate}

As already observed in the introduction this transformation was
firstly introduced in Celestial Mechanics by R. McGehee in its
celebrated paper \cite{McGehee74}. Furthermore this technique became a
milestone in order to investigate the orbital structure of the
$n$-body problem in the neighborhood of total collapse. In fact it was
employed in several problem giving a lot of important feature on the
global dynamics of this interesting problem. However in that context,
due to the fact that the potential is a homogeneous function (of
degree $-1$) allow us to decouple the system in a scalar equation
containing the radial part and in a system which take into account the
angular part. Therefore it makes this technique more feasible for the
application. In fact it is possible to study separately the angular
and radial part; in this perspective it is possible to read the full
behavior of the system by looking only at the angular part. All of
this breakdown in our context since the potential is not homogeneous
anymore.

Another important difference which makes things more involved is that
the equilibrium points are not isolated (in fact they appear in
family) and they do not have a hyperbolic character (and only some of
them are normally hyperbolic manifolds). For instance, to the
knowledge of the authors it is not known if there exists a sort of
Palis Inclination Lemma in this context.

However since the differential of the potential is a homogeneous
function, this still guaranteed the existence of self-similar
(homothetical, homographical) motions.

As direct consequence of the energy relation, for each fixed energy
level $h$ the motion is confined between the two invariant manifold
($\hat E=0$, $r=0$) and the hyperplanes corresponding to the
simultaneously binary collision. However as already observed the
region in which the motion is allowed is not compact. This is
completely different for instance with the case studied by the authors
in \cite{stoica} and in our case a priori unbounded motion can
happen. However as proved above there are no unbounded orbits.

An intriguing question is to understand if this weak logarithmic
singularities always prevents the existence of unbounded motions.



\end{document}